\documentclass[journal,10pt,onecolumn]{IEEEtran}
\usepackage[noadjust]{cite}
\usepackage[pdftex]{graphics}
\usepackage{caption}
\usepackage{subcaption}
\usepackage{amsmath,amsfonts,amssymb}
\usepackage{dsfont}
\usepackage{algorithmic}
\usepackage{array}
\usepackage{stfloats}
\usepackage{bigints}
\usepackage{url}
\usepackage{multirow}
\usepackage{tikz}
\RequirePackage{fix-cm}
\DeclareMathSizes{10}{10}{5}{5}
\newtheorem{theorem}{\textbf{Theorem}}
\newtheorem{proposition}{\textbf{Proposition}}

\begin{document}

\title{Warped Riemannian metrics \\ for location-scale models}

\author{Salem~Said, Lionel Bombrun, Yannick Berthoumieu}

\maketitle

\abstract
The present paper shows that warped Riemannian metrics, a class of Riemannian metrics which play a prominent role in Riemannian geometry, are also of fundamental importance in information geometry. Precisely, the paper features a new theorem, which states that the Rao-Fisher information metric of any location-scale model, defined on a Riemannian manifold, is a warped Riemannian metric, whenever this model is invariant under the action of some Lie group. This theorem is a valuable tool in finding the expression of the Rao-Fisher information metric of location-scale models defined on high-dimensional Riemannian manifolds. Indeed, a warped Riemannian metric is fully determined by only two functions of a single variable, irrespective of the dimension of the underlying Riemannian manifold. Starting from this theorem, several original contributions are made. The expression of the Rao-Fisher information metric of the Riemannian Gaussian model is provided, for the first time in the literature. A generalised definition of the Mahalanobis distance is introduced, which is applicable to any location-scale model defined on a Riemannian manifold. The solution of the geodesic equation is obtained, for any Rao-Fisher information metric defined in terms of warped Riemannian metrics. Finally, using a mixture of analytical and numerical computations, it is shown that the parameter space of the von Mises-Fisher model of $n$-dimensional directional data, when equipped with its Rao-Fisher information metric, becomes a Hadamard manifold, a simply-connected complete Riemannian manifold of negative sectional curvature, for $n = 2,\ldots,8$. Hopefully, in upcoming work, this will be proved for any value of $n$.

\section{Introduction} \label{sec:intro}
\begin{subequations}
Warped Riemannian metrics are a class of Riemannian metrics which arise throughout Riemannian geometry\!~\cite{bishop}\cite{petersen}. For example, the Riemannian metrics of spaces of constant curvature, and of surfaces of revolution, are warped Riemannian metrics. Closely similar to warped Riemannian metrics, warped semi-Riemannian metrics are very important in theoretical physics. Indeed, many gravitational models are given by warped semi-Riemannian metrics\!~\cite{oneil}. The present paper shows that warped metrics, in addition to their well-known role in geometry and physics, play a fundamental role in information geometry, and have a strong potential for applications in statistical inference and statistical learning. 
  
A unified definition of warped Riemannian metrics was first formulated in~\cite{bishop}. Here, only a special case of this definition is required. Precisely, let $M$ be a complete Riemannian manifold, with length element $ds^2_{\scriptscriptstyle M}\,$, and consider the product manifold $\mathcal{M} \,=\, M \times (0\,,\infty)$, equipped with the length element $ds^2_{\scriptscriptstyle \mathcal{M}\,}$
\begin{equation} \label{eq:intro1}
   ds^2_{\scriptscriptstyle\mathcal{M}}(z) \,=\, dr^2 \,+\, \beta^2(r)\,ds^2_{\scriptscriptstyle M}(x) \hspace{1cm} \text{for}\;\; z = (x,r) \in \mathcal{M} 
\end{equation}
where $\beta^2(r)$ is a strictly positive function. Then, the length element $ds^2_{\scriptscriptstyle\mathcal{M}}$ defines a warped Riemannian metric on $\mathcal{M}$. In Riemannian geometry, the coordinate $r$ is a distance function, measuring the distance to some point or hypersurface\!~\cite{petersen}. In physics, $r$ is replaced by the time $t$, and $dr^2$ is replaced by $-dt^2$ in formula (\ref{eq:intro1}) (this is the meaning of ``semi-Riemannian")\!~\cite{oneil}. In any case, the coordinate $x$ can be thought of as a spatial coordinate which determines a position in $M$.  

The intuition behind the present paper is that warped Riemannian metrics are natural candidates for Riemannian metrics on location-scale models. Indeed, if $\mathcal{P}$ is a location-scale model on $M$, with location parameter $\bar{x} \in M$ and scale parameter $\sigma > 0$, then the parameter space of $\mathcal{P}$ is exactly $\mathcal{M} \,=\, M \times (0\,,\infty)$ with its points $z = (\bar{x},\sigma)$. Thus, a warped Riemannian metric on $\mathcal{M}$ can be defined using (\ref{eq:intro1}), after introducing a new scale parameter $r = r(\sigma)$ and setting $x = \bar{x}$. 

As it turns out, this intuition is far from arbitrary. The main new result in the present paper, Theorem \ref{th:warplocdis} of Section \ref{sec:theorem}, states that the Rao-Fisher information metric of any location-scale model is a warped Riemannian metric, whenever this model is invariant under the action of some Lie group. Roughly, Theorem \ref{th:warplocdis} states that if $M$ is a Riemannian symmetric space under the transitive action of a Lie group of isometries $G$, and if each probability density $p(x|\,\bar{x},\sigma)$, belonging to the model $\mathcal{P}$, verifies the invariance condition
\begin{equation} \label{eq:intro2}
  p(\,g\cdot x|\,g\cdot\bar{x}\,,\sigma) \,=\,p(x|\bar{x},\sigma) \hspace{0.25cm} \text{ for all } g \in G
\end{equation}
where $g\cdot x$ denotes the action of $g \in G$ on $x \in M$, then the Rao-Fisher information metric of the model $\mathcal{P}$ is a warped Riemannian metric. A technical requirement for Theorem \ref{th:warplocdis} is that the Riemannian symmetric space $M$ should be irreducible. The meaning of this requirement, and the fact that it can be relaxed in certain cases, are discussed in Remarks 4 and 5 of Section \ref{sec:theorem}. The proof of Theorem \ref{th:warplocdis} is given in Appendix \ref{app:A}. 

A fundamental idea of information geometry is that the parameter space of a statistical model $\mathcal{P}$ should be considered as a Riemannian manifold\!~\cite{amari}\cite{chentsov}. According to\!~\cite{chentsov}, the unique way of doing so is by turning Fisher's information matrix into a Riemannian metric, the Rao-Fisher information metric. In this connection, Theorem \ref{th:warplocdis} shows that, when the statistical model $\mathcal{P}$ is a location-scale model which is invariant under the action of a Lie group, information geometry inevitably leads to the study of warped Riemannian metrics. 

In addition to stating and proving Theorem \ref{th:warplocdis}, the present paper aims to explore its implications, with regard to the Riemannian geometry of location-scale models, and to lay the foundation for its applications in statistical inference and statistical learning.

To begin, Section \ref{sec:examples} applies Theorem \ref{th:warplocdis} to two location-scale models, the von Mises-Fisher model of directional data~\!\cite{jupp}\cite{chikuse}, and the Riemannian Gaussian model of data in spaces of covariance matrices~\!\cite{vemuri,said1,said2}. This leads to the analytic expression of the Rao-Fisher information metric of each one of these two models. Precisely, the Rao-Fisher information metric of the von Mises-Fisher model is given in Proposition \ref{prop:vmf}, and that of the Riemannian Gaussian model is given in Proposition \ref{prop:rgm}. The result of Proposition \ref{prop:vmf} is essentially already contained in~\cite{jupp}, (see Page 199), but Proposition \ref{prop:rgm} is new in the literature. 

Finding the analytic expression of the Rao-Fisher information metric, or equivalently of Fisher's information matrix, of a location-scale model $\mathcal{P}$ defined on a high-dimensional non-trivial manifold $M$, is a very difficult task when attempted by direct calculation. Propositions \ref{prop:vmf} and \ref{prop:rgm} show that this task is greatly simplified by Theorem \ref{th:warplocdis}. Precisely, if the dimension of $M$ is $d$, then the dimension of the parameter space $\mathcal{M} \,=\, M \times (0\,,\infty)$ is $d+1$. Therefore, \textit{a priori}, the expression of the Rao-Fisher information metric involves $(d+1)(d+2)/2$ functions of both parameters $\bar{x}$ and $\sigma$ of the model $\mathcal{P}$. Instead of so many functions of both $\bar{x}$ and $\sigma$, Theorem \ref{th:warplocdis} reduces the expression of the Rao-Fisher information metric to only two functions of $\sigma$ alone. In the notation of (\ref{eq:intro1}), these two functions are $\alpha(\sigma) = dr/d\sigma$ and $\beta(\sigma) = \beta(r(\sigma))$. 

Section \ref{sec:maha} builds on Theorem \ref{th:warplocdis} to introduce a general definition of the Mahalanobis distance, applicable to any location-scale model $\mathcal{P}$ defined on a manifold $M$. Precisely, assume that the model $\mathcal{P}$ verifies the conditions of Theorem \ref{th:warplocdis}, so its Rao-Fisher information metric is a warped Riemannian metric. Then, the \textit{generalised Mahalanobis distance} is defined as the Riemannian distance on $M$ which is induced by the restriction of the Rao-Fisher information metric to $M$. The expression of the generalised Mahalanobis distance is given in Propositions \ref{prop:maha1} and \ref{prop:maha2}. It was recently applied to visual content classification in~\cite{icip}. 

The generalised Mahalanobis distance includes the classical Mahalanobis distance as a special case. Precisely, assume $\mathcal{P}$ is the isotropic normal model defined on $M = \mathbb{R}^d$, so each density $p(x|\,\bar{x},\sigma)$ is a $d$-variate normal density with mean $\bar{x}$ and covariance matrix $\sigma^2$ times the identity. Then, $\mathcal{P}$ verifies the invariance condition (\ref{eq:intro2}) under the action of the group $G$ of translations in $\mathbb{R}^d$. Therefore, by Theorem \ref{th:warplocdis}, its Rao-Fisher information metric is a warped Riemannian metric. This metric is already known in the literature, in terms of the length element\!~\cite{bensadon}\cite{bensadonth}
\begin{equation} \label{eq:poincare1}
 ds^2_{\scriptscriptstyle \mathcal{M}}(z) = \frac{2d}{\sigma^2}\,d\sigma^2\,+\,\frac{1}{\sigma^2}\,\Vert d\bar{x}\Vert^2 \hspace{1cm} \text{for}\;\; z = (\bar{x},\sigma) \in \mathcal{M} 
\end{equation}
where $\Vert d\bar{x}\Vert^2$ denotes the Euclidean length element on $\mathbb{R}^d$. Now, the restriction of this Rao-Fisher information metric to $M = \mathbb{R}^d$ is given by the second term in (\ref{eq:poincare1}), which is clearly the Euclidean length element divided by $\sigma^2$, and corresponds to the classical Mahalanobis distance\!~\cite{maclahlan} --- note that (\ref{eq:poincare1}) is brought into the form (\ref{eq:intro1}) by letting $r(\sigma) = {\scriptstyle \sqrt{2d}}\,\log(\sigma)$ and $\beta(r) = \exp(-r/{\scriptstyle \sqrt{2d}})$.

Section \ref{sec:pn} illustrates the results of Sections \ref{sec:examples} and \ref{sec:maha}, by applying them to the special case of the Riemannian Gaussian model defined on $M = \mathcal{P}_n\,$, the space of $n\times n$ real covariance matrices. In particular, it gives directly applicable expressions of the Rao-Fisher information metric, and of the generalised Mahalanobis distance, corresponding to this model. As it turns out, the generalised Mahalanobis distance defines a whole new family of affine-invariant distances on $\mathcal{P}_n\,$, in addition to the usual affine-invariant distance, which was introduced to the information science community in\!~\cite{pennec2}. 

Section \ref{sec:geodesic} provides the solution of the geodesic equation of any of the Rao-Fisher information metrics arising from Theorem \ref{th:warplocdis}. The main result is Proposition \ref{prop:geodesicmultwarp}, which states that the solution of this equation, for given initial conditions, reduces to the solution of a one-dimensional second order differential equation. This implies that geodesics with given initial conditions can be constructed at a reasonable numerical cost, which opens the possibility, with regard to future work, of a practical implementation of Riemannian line-search optimisation algorithms, which find the extrema of cost functions by searching for them along geodesics\!~\cite{absil}\cite{lue}.
 
Section \ref{sec:hadamard} is motivated by the special case of the isotropic normal model, with its Rao-Fisher information metric given by (\ref{eq:poincare1}). It is well-known that, after a trivial change of coordinates, the length element (\ref{eq:poincare1}) coincides with the length element of the Poincar\'e half-space model of hyperbolic geometry\!~\cite{bensadon}\cite{bensadonth}. This means that the parameter space of the isotropic normal model, when equipped with its Rao-Fisher information metric, becomes a space of constant negative curvature, and in particular a Hadamard manifold, a simply-connected complete Riemannian manifold of negative sectional curvature\!~\cite{petersen}\cite{chavel}. One cannot but wonder whether other location-scale models also give rise to Hadamard manifolds in this way. This is investigated using Propositions \ref{prop:vmf} and \ref{prop:geodesicmultwarp}, for the case of the von Mises-Fisher model. A mixture of analytical and numerical computations leads to a surprising new observation\,: the parameter space of the von Mises-Fisher model of $n$-dimensional directional data, when equipped with its Rao-Fisher information metric, becomes a Hadamard manifold, for $n = 2,\ldots,8$. Ongoing research warrants the conjecture that this is true for any value of $n$, but this is yet to be fully proved. 

Theorem \ref{th:warplocdis}, the main new result in the present paper, has many potential applications, which will be developed in future work. Indeed, this theorem can provide the expression of the Rao-Fisher information metric, or equivalently of Fisher's information matrix, of a location-scale model, even if this model is defined on a high-dimensional non-trivial manifold. By doing so, it unlocks access to the many applications which require this expression, both in statistical inference and in statistical learning. In statistical inference, the expression of the Rao-Fisher information metric allows the computation of the Cram\'er-Rao lower bound, and the construction of asymptotic chi-squared statistics\!~\cite{lehman}\cite{inference}. In statistical learning, it allows the practical implementation of the natural gradient algorithm, which has the advantages of full reparameterisation invariance, of asymptotic efficiency of its stochastic version, and of the linear rate of convergence of its deterministic version~\cite{bottou,amaring,martens}.

A first step, towards developing the applications of Theorem \ref{th:warplocdis}, was recently taken in~\cite{paolomix}. Using the expression of the Rao-Fisher information metric of the Riemannian Gaussian model, this derived and successfully implemented the natural gradient algorithm, for the problem of on-line learning of an unknown probability density, on the space $\mathcal{P}_n$ of $n\times n$ real covariance matrices. Future work will focus on extending this algorithm to more problems of statistical learning, including problems of on-line classification and regression in the  space $\mathcal{P}_{n\,}$, and in other spaces of covariance matrices. In addition to its use in deriving the natural gradient algorithm for problems of statistical learning, the expression of the Rao-Fisher information metric of the Riemannian Gaussian model can be used in deriving so-called natural evolutionary strategies, for black-box optimisation in spaces of covariance matrices.  These would generalize currently existing natural evolutionary strategies, which are mostly restricted to black-box optimisation in Euclidean space\!~\cite{nes}\cite{igo}. 

The following Section \ref{sec:warped} provides necessary background on warped Riemannian metrics. The reader will readily notice the difference in notation, between the present introduction and Section \ref{sec:warped}. Here, in formulae (\ref{eq:intro1}) and (\ref{eq:poincare1}), Riemannian metrics were given in terms of length elements. In Section \ref{sec:warped}, Riemannian metrics will be given by scalar products on the tangent space. It was felt that the use of length elements, while somewhat less rigorous, is more intuitive, and therefore better suited for the introduction.
 
\end{subequations}

\vfill
\pagebreak

\section{Background on warped Riemannian metrics} \label{sec:warped}
\begin{subequations} \label{eq:def}
Assume $M$ is a complete Riemannian manifold with Riemannian metric $Q$, and consider the manifold $\mathcal{M} \,= M\times (0\,,\infty)$. A warped Riemannian metric $I$ on $\mathcal{M}$ is given in the following way~\cite{bishop,oneil,petersen}. Let $\alpha$ and $\beta$ be positive functions, defined on $(0\,,\infty)$. Then, for $z = (x,\sigma) \in \mathcal{M}$, let the scalar product $I_z$ on the tangent space $T_z\mathcal{M}$ be defined by
\begin{eqnarray} \label{eq:warpdef}
 I_z(U,U) \,=\, \left(\,\!\alpha(\sigma) \, u_\sigma\,\!\right)^2 \,+\, \beta^2(\sigma)\,Q_x(u,u)  \hspace{1cm} U \in T_z\mathcal{M}
\end{eqnarray}
where $U \,=\, u_\sigma\,\partial_\sigma\,+\,u$ with $u_\sigma \in \mathbb{R}$ and $u \in T_xM$.  The functions $\alpha$ and $\beta$ are part of the definition of the warped metric $I$. Once these functions are fixed, it is possible to introduce a change of coordinates $r = r(\sigma)$ which eliminates $\alpha$ from (\ref{eq:warpdef}). Precisely, if $dr/d\sigma=\,\alpha(\sigma)$ then 
\begin{eqnarray} \label{eq:warpdefr}
 I_z(U,U) \,=\, u^2_r \,+\, \beta^2(r)\,Q_x(u,u)  
\end{eqnarray}
where $U \,=\, u_r\,\partial_r\,+\,u$ and $\beta(r) = \beta(\sigma(r))$. 
\end{subequations}

\begin{subequations} \label{eq:r}
The coordinate $r$ will be called \textit{vertical distance}. This is not a standard terminology, but is suggested as part of the following  geometric picture. For $z = (x,\sigma) \in \mathcal{M}$, think of $x$  as a horizontal coordinate, and of $\sigma$ as a vertical coordinate. Accordingly, the points $z_{\scriptscriptstyle 0} = (x,\sigma_{\scriptscriptstyle 0})$ and $z_{\scriptscriptstyle 1} = (x,\sigma_{\scriptscriptstyle 1})$ lie on the same vertical line. It can be shown from (\ref{eq:warpdefr}) that the Riemannian distance between $z_{\scriptscriptstyle 0}$ and $z_{\scriptscriptstyle 1}$ is
\begin{equation} \label{eq:verticald}
 d\,\!(z_{\scriptscriptstyle 0}\,\!,z_{\scriptscriptstyle 1}) \,=\, r(\sigma_{\scriptscriptstyle 1}) - r(\sigma_{\scriptscriptstyle 0}) \hspace{1cm} \text{where } \,\,\sigma_{\scriptscriptstyle 0\,} <\sigma_{\scriptscriptstyle 1\,}
\end{equation}
Precisely, $d\,\!(z_{\scriptscriptstyle 0}\,\!,z_{\scriptscriptstyle 1})$ is the Riemannian distance induced by the warped Riemannian metric $I$. 

The vertical distance $r$ can be used to express a necessary and sufficient condition for completeness of the manifold $\mathcal{M}$, equipped with the warped Riemannian metric $I$. Namely, $\mathcal{M}$ is a complete Riemannian manifold, if and only if
\begin{equation} \label{eq:complete}
 \lim_{\sigma\rightarrow \infty} r(\sigma)- r(\sigma_{\scriptscriptstyle 0}) = \infty \;\;\text{and}\;\;
\lim_{\sigma\rightarrow 0} r(\sigma_{\scriptscriptstyle 1})- r(\sigma) = \infty
\end{equation} 
where $\sigma_{\scriptscriptstyle 0}$ and $\sigma_{\scriptscriptstyle 1}$ are arbitrary. This condition is a special case of Lemma 7.2 in~\cite{bishop}.
\end{subequations}

Let $K^\mathcal{M}$ and $K^M$ denote the sectional curvatures of $\mathcal{M}$ and $M$, respectively. The relation between these two is given by the curvature equations of Riemannian geometry~\cite{petersen}\cite{docarmo}. These are,   
\begin{subequations} \label{eq:curvature}  
\begin{eqnarray} 
\label{eq:gauss} \text{Gauss equation\,:}  &\hspace{0.25cm} K^\mathcal{M}_z(u,v) \,= \frac{1}{\strut\beta^2}\,K^M_x(u,v) -  \left(\frac{\partial_r\beta}{\strut\beta}\right)^2 & \hspace{0.2cm} u,v \in T_xM \\[0.1cm]
\label{eq:radial} \text{radial curvature equation\,:} &\hspace{0.25cm} K^\mathcal{M}_z(u,\partial_r) \,= - \frac{\partial^2_r\beta}{\strut\beta} \hspace{2.3cm} 
\end{eqnarray}
 \end{subequations}
Here, the notations $K^\mathcal{M}_z$ and $K^M_x$ mean that $K^\mathcal{M}$ is computed at $z$, and $K^M$ is computed at $x$, where $z = (x,\sigma)$. Equations (\ref{eq:curvature}) are a special case of Lemma 7.4 in~\cite{bishop}. 

Note, as a corollary of these equations, that $\mathcal{M}$ has negative sectional curvature $K^\mathcal{M} < 0$, if $M$ has negative sectional curvature $K^M < 0$ and $\beta$ is a strictly convex function of $r$. \\[0.1cm]
\textbf{Remark 1\,:} equations (\ref{eq:def}) contains an abuse of notation. Namely, $u$ denotes a tangent vector to $M$ at $x$, and a tangent vector to $\mathcal{M}$ at $z$, at the same time. In the mathematical literature (for example, in~\cite{bishop}\cite{oneil}), one writes $d\pi_z(U)$ instead of $u$, using the derivative $d\pi$ of the projection mapping $\pi(z) = x$, and this eliminates any ambiguity. In the present paper, a deliberate choice is made to use a lighter, though not entirely correct, notation. \hfill$\blacksquare$\\[0.1cm]
\textbf{Remark 2\,:} consider the proof of equations (\ref{eq:r}). For (\ref{eq:verticald}), let $\gamma(t)$ and $c(t)$ be curves connecting $z_{\scriptscriptstyle 0} = (x,\sigma_{\scriptscriptstyle 0})$ and $z_{\scriptscriptstyle 1} = (x,\sigma_{\scriptscriptstyle 1})$. Assume these are parameterised as follows, 
$$
\gamma(t)\,:\,\left\lbrace\!\! \begin{array}{ll}
x(t) \hspace{-0.1cm}& \!\!= x \, \, \text{ (constant)}  \\[0.1cm]
r(t) \hspace{-0.3cm}& \!\!= r(\sigma_{\scriptscriptstyle 0}) \,+\, t\, (r(\sigma_{\scriptscriptstyle 1})-r(\sigma_{\scriptscriptstyle 0}))
\end{array}\!\!\right\rbrace
\hspace{1.2cm}
c(t)\,:\,\left\lbrace\!\! \begin{array}{l}
x(t) \\[0.1cm]
r(t)
\end{array}\!\!\right\rbrace
$$
where $t \in [0,1]$. If $L(\gamma)$ and $L(c)$ denote the lengths of these curves, then from (\ref{eq:warpdefr}),
$$
L(c) = \int^1_0\, \left( \left(\dot{r}\right)^2 \,+\, \beta^2(r) Q(\dot{x},\dot{x})\right)^{1/2}\,dt \,\geq\, \int^1_0\, \dot{r}\,dt
\,=\, r(\sigma_{\scriptscriptstyle 1})-r(\sigma_{\scriptscriptstyle 0}) \,=\, L(\gamma)$$
where the dot denotes differentiation with respect to $t$, and the inequality is strict unless $\gamma = c$. This shows that $\gamma(t)$ is the unique 
length-minimising geodesic connecting $z_{\scriptscriptstyle 0}$ and $z_{\scriptscriptstyle 1}$. Thus, $d\,\!(z_{\scriptscriptstyle 0}\,\!,z_{\scriptscriptstyle 1}) = L(\gamma)$, and this gives (\ref{eq:verticald}). For (\ref{eq:complete}), note that Lemma 7.2 in~\cite{bishop} states that $\mathcal{M}$ is complete, if and only if $(0\,,\infty)$ is complete when equipped with the distance
$$
d_{(0,\infty)}(\sigma_{\scriptscriptstyle 0}\,,\sigma_{\scriptscriptstyle 1}) \,= \left| r(\sigma_{\scriptscriptstyle 1}) - r(\sigma_{\scriptscriptstyle 0})\right|
$$
However, this is equivalent to (\ref{eq:complete}). An alternative way of understanding (\ref{eq:complete}) is provided in Remark 14 of Section \ref{sec:geodesic}. Precisely, let a geodesic of the form $\gamma(t)$ be called a vertical geodesic. The first condition in (\ref{eq:complete}) states that a vertical geodesic cannot reach the value $\sigma = \infty$ within a finite time, and the second condition in (\ref{eq:complete}) states that a vertical geodesic cannot reach the value $\sigma = 0$ within a finite time. \hfill$\blacksquare$ \\[0.1cm] 
\textbf{Remark 3\,:} for each $\sigma \in (0\,,\infty)$, the manifold $M$ can be embedded into the manifold $\mathcal{M}$, in the form of the hypersurface $M_\sigma = M \times \lbrace \sigma \rbrace$. Through this embedding, the warped Riemannian metric $I$ of $\mathcal{M}$ induces a Riemannian metric $Q^\sigma$ on $M$. By definition, this metric $Q^\sigma$ is obtained by the restriction of $I$ to the tangent vectors of $M_\sigma$\!~\cite{petersen}\cite{docarmo}. It follows from (\ref{eq:def}) that
\begin{equation} \label{eq:inducedmetric}
  Q^\sigma_x(u,u) \,=\, \beta^2(\sigma)\,Q_x(u,u)  
\end{equation}
The induced metric $Q^\sigma$ will be called an extrinsic metric on $M$, since it comes from the ambient space $\mathcal{M}$. By (\ref{eq:inducedmetric}),
the extrinsic metric $Q^\sigma$ is equal to a scaled version of the Riemannian metric $Q$ of $M$, with scaling factor equal to $\beta(\sigma)$. \hfill$\blacksquare$
\section{Connection with location-scale models} \label{sec:theorem}
 This section establishes the connection between warped Riemannian metrics and location-scale models. The main result is Theorem \ref{th:warplocdis}, which states that the Rao-Fisher information metric of any location-scale model is a warped Riemannian metric, whenever this model is invariant under the action of some Lie group.

To state this theorem, assume $M$ is an irreducible Riemannian symmetric space, with invariant Riemannian metric $Q$, under the transitive action of a Lie group of isometries $G$~\cite{helgason}. Consider a location-scale model $\mathcal{P}$ defined on $M$,
\begin{equation} \label{eq:locationscale}
  \mathcal{P} = \left\lbrace\,p(x|\bar{x},\sigma)\,;\,\bar{x} \in M \,,\, \sigma \in (0\,,\infty)\right\rbrace
\end{equation}
To each point $z = (\bar{x},\sigma)$ in the parameter space $\mathcal{M} = M \times (0\,,\infty)$, this model associates a probability density $p(x|\bar{x},\sigma)$ on $M$, which has a location parameter $\bar{x}$ and a scale parameter $\sigma$. Precisely, $p(x|\bar{x},\sigma)$ is a probability density with respect to the invariant Riemannian volume element of $M$. 

The condition that the model $\mathcal{P}$ is invariant under the action of the Lie group $G$ means that, 
\begin{equation} \label{eq:invariance}
  p(\,g\cdot x|\,g\cdot\bar{x}\,,\sigma) \,=\,p(x|\bar{x},\sigma) \hspace{0.25cm} \text{ for all } g \in G
\end{equation}
where $g\cdot x$ denotes the action of $g \in G$ on $x \in M$. 

The Rao-Fisher information metric of the location-scale model $\mathcal{P}$ is a Riemannian metric $I$ on the parameter space $\mathcal{M}$ of this model~\cite{amari}.  It is defined as follows, for $z =(\bar{x},\sigma) \in \mathcal{M}$ and $U\in T_z\mathcal{M}$,
\begin{equation} \label{eq:rao}
  I_z(U,U) \,=\, \mathbb{E}_z\left( \,d\ell(z)\,U\,\right)^2 
\end{equation}
where $\mathbb{E}_z$ denotes expectation with respect to the probability density $p(x|z) = p(x|\bar{x},\sigma)$, and $\ell(z)$ is the log-likelihood function,
given by $\ell(z)(x) = \log p(x|z)$.  

In the following statement, $\nabla_{\bar{x}}\,\ell(z)$ denotes the Riemannian gradient of $\ell(z)$, taken with respect to $\bar{x} \in M$, while the value of $\sigma$ is fixed. 
\begin{theorem} \label{th:warplocdis}
if condition (\ref{eq:invariance}) is verified, then the Rao-Fisher information metric $I$ of (\ref{eq:rao}) is a warped Riemannian metric given by (\ref{eq:warpdef}), where  
\begin{equation} \label{eq:warplocdis}
  \alpha^2(\sigma) \,=\mathbb{E}_z\left(\partial_\sigma\ell(z)\right)^2 \hspace{1cm} \beta^2(\sigma) \,= \left. \mathbb{E}_z\, Q\left(\nabla_{\bar{x}}\,\ell(z)\,,\nabla_{\bar{x}}\,\ell(z)\,\right) \middle/ \mathrm{dim}\,M \right.
\end{equation}
The expectations appearing in (\ref{eq:warplocdis}) do not depend on $\bar{x}$, so $\alpha(\sigma)$ and $\beta(\sigma)$ are well-defined functions of $\sigma$.
\end{theorem}
\vspace{0.13cm}
\textbf{Remark 4\,:} recall the definition of an irreducible Riemannian symmetric space~\cite{helgason}. A Riemannian manifold $M$, whose group of isometries is denoted $G$, is called a Riemannian symmetric space, if for each $\bar{x} \in M$ there exists an isometry $s_{\bar{x}} \in G$, whose effect is to fix $\bar{x}$ and reverse the geodesic curves passing through $\bar{x}$. Further, $M$ is called irreducible if it verifies the following condition. Let $K_{\bar{x}}$ be the subgroup of $G$ which consists of those elements $k$ such that $k\cdot\bar{x} = \bar{x}$. For each $k \in K_{\bar{x}}$, its derivative $dk_{\bar{x}}$ is a linear mapping of $T_{\bar{x}}M$. The mapping $k \mapsto dk_{\bar{x}}$ is a representation of $K_{\bar{x}}$ in $T_{\bar{x}}M$, called the isotropy representation, and $M$ is called an irreducible Riemannian symmetric space if the isotropy representation is irreducible. That is, if the isotropy representation has no invariant subspaces in $T_{\bar{x}}M$, except $\lbrace 0 \rbrace$ and $T_{\bar{x}}M$. Irreducible Riemannian symmetric spaces are classified in~\cite{helgason} (Table I, Page 346 and Table II, Page 354). They include spaces of constant curvature, such as spheres and hyperbolic spaces, as well as spaces of positive definite matrices which have determinant equal to $1$, and whose entries are real or complex numbers, or quaternions. \hfill$\blacksquare$\\[0.1cm]
\textbf{Remark 5\,:} it is sometimes possible to apply Theorem \ref{th:warplocdis}, even when the Riemannian symmetric space $M$ is not irreducible. For example, in \ref{subsec:gauss}, Theorem \ref{th:warplocdis} will be used to find the expression of the Rao-Fisher information metric of the Riemannian Gaussian model~\cite{vemuri,said1,said2}. For this model, when $M$ is not irreducible, the Rao-Fisher information metric turns out to be a so-called multiply-warped Riemannian metric, rather than a warped Riemannian metric. The concrete case of $M = \mathcal{P}_n$, the space of $n \times n$ real covariance matrices, is detailed in Section \ref{sec:pn}.\hfill$\blacksquare$\\[0.1cm]
\textbf{Proof of Theorem \ref{th:warplocdis}\,:} 
recall the expression $U \,=\, u_\sigma\,\partial_\sigma\,+\,u$ with $u_\sigma \in \mathbb{R}$ and $u \in T_{\bar{x}}M$. Since the Rao-Fisher information metric $I$ is bilinear and symmetric,
$$
I_z(U,U) \,=\, I_z(\partial_\sigma,\partial_\sigma)\,u^2_\sigma \,+\, 2I_z(\partial_\sigma,u)\,u_\sigma\,+\,I_z(u,u)
$$
It is possible to show the following,
\begin{subequations} \label{eq:proof}
\begin{eqnarray}
\label{eq:proof1} I_z(\partial_\sigma,\partial_\sigma) =& \alpha^2(\sigma) \\[0.12cm]
\label{eq:proof2} I_z(\partial_\sigma,u) =&  0 \\[0.12cm]
\label{eq:proof3} I_z(u,u) =& \beta^2(\sigma)\,Q_{\bar{x}}(u,u)  
\end{eqnarray}
\end{subequations}
where $\alpha^2(\sigma)$ and $\beta^2(\sigma)$ are given by (\ref{eq:warplocdis}). \\[0.1cm]
\textit{Proof of (\ref{eq:proof1})\,:} this is immediate from (\ref{eq:rao}). Indeed, 
$$
I_z(\partial_\sigma,\partial_\sigma) = \mathbb{E}_z\left( \,d\ell(z)\,\partial_\sigma\,\right)^2 = \mathbb{E}_z\left(\partial_\sigma\ell(z)\right)^2 
$$
\textit{Proof of (\ref{eq:proof2})\,:} this is carried out in Appendix \ref{app:A}, using the fact that $M$ is a Riemannian symmetric space.  \\[0.1cm]
\textit{Proof of (\ref{eq:proof3})\,:} this is carried out in Appendix \ref{app:A}, using the fact that $M$ is irreducible, by an application of Schur's lemma from the theory of group representations~\cite{knapp}.

The fact that the expectations appearing in (\ref{eq:warplocdis}) do not depend on $\bar{x}$ is also proved in Appendix \ref{app:A}. Throughout the proof of the theorem, the following identity is used, which is equivalent to condition (\ref{eq:invariance}). For any real-valued function $f$ on $M$,
\begin{equation} \label{eq:fcirc}
\mathbb{E}_{g\cdot z} \,f =  \mathbb{E}_z\left(f\circ g\right)
\end{equation}
Here, $g\cdot z = (\,g\cdot \bar{x},\sigma)$, and $f \circ g$ is the function $(f\circ g)(x) = f(g\cdot x)$, for $g \in G$ and $z = (\bar{x},\sigma)$.  \hfill$\blacksquare$



%



\section{Examples\,: von Mises-Fisher and Riemannian Gaussian} \label{sec:examples}
This section applies Theorem \ref{th:warplocdis} to finding the expression of the Rao-Fisher information metric of two location-scale models. These are the von Mises-Fisher model, which is widely used in the study of directional data~\cite{jupp}\cite{chikuse}, and the Riemannian Gaussian model, recently introduced in the study of data with values in spaces of covariance matrices~\cite{vemuri,said1,said2}. 

The application of Theorem \ref{th:warplocdis} to these two models is encapsulated in the following Proposition \ref{prop:exponential}. Precisely, both of these models are of a common exponential form, which can be described as follows. Let $M$ be an irreducible Riemannian symmetric space, as in Section \ref{sec:theorem}. In the notation of (\ref{eq:locationscale}), consider a location-scale model $\mathcal{P}$ defined on $M$, by
\begin{subequations}\label{eq:exponential}
\begin{equation} \label{eq:exponential1}
  p(x|\bar{x},\sigma) \,=\, \exp\left[\,\eta(\sigma)\,D(x,\bar{x}) \,-\, \psi(\eta(\sigma))\right]
\end{equation}
where $\eta(\sigma)$ is a certain parameter, to be called the natural parameter, and where $D: M \times M \rightarrow \mathbb{R}$ verifies the condition, 
\begin{equation} \label{eq:exponential2}
  D(\,g\cdot x\,,\,\,g\cdot\bar{x}) \,=\,D(x\,,\,\bar{x}) \hspace{0.25cm} \text{ for all } g \in G
\end{equation}
There is no need to assume that the function $D$ is positive.
\end{subequations}
\begin{proposition} \label{prop:exponential}
  if the model $\mathcal{P}$ is given by equations (\ref{eq:exponential}), then the Rao-Fisher information metric $I$ of this model is a warped Riemannian metric,
\begin{subequations} \label{eq:propexp} 
\begin{equation} \label{eq:propexp1} 
   I_z(U,U) \,=\, \psi^{\prime\prime}(\eta)\,u^2_\eta\,+\, \beta^2(\eta)\,Q_{\bar{x}}(u,u) 
\end{equation}
where $U \,=\, u_\eta\,\partial_\eta\,+\,u$ with $u_\eta \in \mathbb{R}$ and $u \in T_{\bar{x}}M$, and where
\begin{equation} \label{eq:propexp2}
 \beta^2(\eta) \,= \,\eta^2\left. \mathbb{E}_z\,Q\left(\nabla_{\bar{x}}\,D,\nabla_{\bar{x}}\,D\,\right) \middle/ \mathrm{dim}\,M \right.
\end{equation}
\end{subequations}
\end{proposition}
\textbf{Proof\,:} for a model $\mathcal{P}$ defined by (\ref{eq:exponential1}), condition (\ref{eq:exponential2}) is equivalent to condition (\ref{eq:invariance}). Therefore, by application of Theorem \ref{th:warplocdis}, it follows that $I$ is a warped Riemannian metric, of the form (\ref{eq:warpdef}),
\begin{subequations} \label{eq:proff}
\begin{equation} \label{eq:proff1}
I_z(U,U) \,=\, \left(\,\!\alpha(\sigma) \, u_\sigma\,\!\right)^2 \,+\, \beta^2(\sigma)\,Q_{\bar{x}}(u,u)
\end{equation}
where $\alpha^2(\sigma)$ and $\beta^2(\sigma)$ are given by (\ref{eq:warplocdis}). Consider the first term in (\ref{eq:proff1}). By the change of coordinates formula~\cite{lee}, $u_\sigma = \sigma^\prime(\eta)\,u_\eta$, where the prime denotes differentiation with respect to $\eta$. It follows that
\begin{equation} \label{eq:proff2}
\left(\,\!\alpha(\sigma) \, u_\sigma\,\!\right)^2 \,=\, \alpha^2(\sigma) \left(\sigma^\prime(\eta)\right)^{\!2}\, u^2_\eta
\end{equation}
However, by (\ref{eq:warplocdis}),
\begin{equation} \label{eq:proff3}
\alpha^2(\sigma) \left(\sigma^\prime(\eta)\right)^2 \,=\,  \mathbb{E}_z\left(\, \partial_\sigma\ell(z)\,\sigma^\prime(\eta)\, \right)^2 \,=\, 
\mathbb{E}_z\left(\, \partial_\eta\ell(z)\, \right)^2
\end{equation}
Here, the log-likelihood $\ell(z)$ is found from (\ref{eq:exponential1}),
\begin{equation} \label{eq:expll}
 \ell(z)(x) = \eta(\sigma)\,D(x,\bar{x}) \,-\, \psi(\eta(\sigma))
\end{equation}
Therefore, the last expression in (\ref{eq:proff3}) is
\begin{equation} \label{eq:amarari}
\mathbb{E}_z\left(\, \partial_\eta\ell(z)\, \right)^2 \,= \, - \mathbb{E}_z\, \partial^2_\eta \ell(z) \,=\, \psi^{\prime\prime}(\eta)
\end{equation}
where the first equality is the same as in~\cite{amari}, (see Page 28). Now, (\ref{eq:proff2}) and (\ref{eq:proff3}) imply
\begin{equation} \label{eq:proff4}
 \left(\,\!\alpha(\sigma) \, u_\sigma\,\!\right)^2 \,=\, \psi^{\prime\prime}(\eta)\,u^2_\eta
\end{equation}
Replacing this in (\ref{eq:proff1}), and writing $\beta(\eta) = \beta(\sigma(\eta))$, gives
\begin{equation} \label{eq:proff5}
I_z(U,U) \,=\, \psi^{\prime\prime}(\eta)\,u^2_\eta \,+\, \beta^2(\eta)\,Q_{\bar{x}}(u,u)
\end{equation}
which is the same as (\ref{eq:propexp1}). To prove the proposition, it remains to show that $\beta^2(\eta)$ is given by (\ref{eq:propexp2}). To do so, note that it follows from (\ref{eq:expll}),
$$
\nabla_{\bar{x}}\,\ell(z) = \nabla_{\bar{x}}\left[\,\eta(\sigma)\,D(x,\bar{x}) \,-\, \psi(\eta(\sigma))\right] = 
\eta(\sigma)\,\nabla_{\bar{x}}\,D(x,\bar{x}) 
$$
Replacing this in (\ref{eq:warplocdis}) gives,
$$
\beta^2(\eta) \,= \left. \mathbb{E}_z\, Q\left(\nabla_{\bar{x}}\,\ell(z)\,,\nabla_{\bar{x}}\,\ell(z)\,\right) \middle/ \mathrm{dim}\,M \right. \,= 
\eta^2\, \left. \mathbb{E}_z\,Q\left(\nabla_{\bar{x}}\,D,\nabla_{\bar{x}}\,D\,\right) \middle/ \mathrm{dim}\,M \right.
$$
\end{subequations}
and this is the same as (\ref{eq:propexp2}). \hfill$\blacksquare$
\subsection{The von Mises-Fisher model} \label{subsec:vmf}
The von Mises-Fisher model is a mainstay of directional statistics~\cite{jupp}\cite{chikuse}. In the notation of  (\ref{eq:exponential}), this model corresponds to $M = S^{n-1}$, the unit sphere in $\mathbb{R}^{n}$, and to $G = O(n)$, the group of $n\times n$ real orthogonal matrices, which acts on $\mathbb{R}^n$ by rotations. Then, the expressions appearing in (\ref{eq:exponential1}) are
\begin{equation} \label{eq:vmf}
   D(x,\bar{x}) = \langle x,\bar{x}\rangle \hspace{1cm} \psi(\eta) = \nu\log(2\pi) \,+\,\log\left( \eta^{1-\nu}I_{\nu-1}(\eta)\right) \text{ for } \eta \in [\,0\,,\infty)
\end{equation}
where $\langle x,\bar{x}\rangle$ denotes the Euclidean scalar product in $\mathbb{R}^n$, so that condition (\ref{eq:exponential2}) is clearly verified, and where $\nu = n/2$ and $I_{\nu - 1}$ denotes the modified Bessel function of order $\nu - 1$. Here, the natural parameter $\eta$ and the scale parameter $\sigma$ should be considered identical, in the sense that $\eta(\sigma) = \sigma$, as long as $\sigma \in (0\,,\infty)$. However, $\eta$ takes on the additional value $\eta = 0$, which requires a special treatment.\\[0.1cm]
\textbf{Remark 6\,:} the parameter space of the von Mises-Fisher model will be identified with the space $\mathbb{R}^n$. This is done by mapping each couple $(\bar{x},\eta)$ to the point $z =\eta\,\bar{x}$ in $\mathbb{R}^n$. This mapping defines a diffeomorphism from the set of couples $(\bar{x},\eta)$ where $\eta \in (0\,,\infty)$, to the open subset $\mathbb{R}^n - \lbrace 0 \rbrace \,\subset \mathbb{R}^n$. On the other hand, it maps all couples 
$(\bar{x},\eta=0)$, to the same point $z = 0 \in \mathbb{R}^n$. Note that each couple $(\bar{x},\eta)$ where $\eta \in (0\,,\infty)$ defines a distinct von Mises-Fisher distribution, which is a unimodal distribution with its mode at $\bar{x}$. On the other hand, all couples $(\bar{x},\eta=0)$ define the same von Mises-Fisher distribution, which is the uniform distribution on $S^{n-1}$. Therefore, it is correct to map all of these couples to the same point $z = 0$. \hfill$\blacksquare$
\vfill
\pagebreak
Proposition \ref{prop:exponential} will only provide the Rao-Fisher information metric of the von Mises-Fisher model on the subset $\mathbb{R}^n - \lbrace 0 \rbrace$ of the parameter space $\mathbb{R}^n$. Therefore, it is necessary to verify that this metric has a well-defined limit at the point $z = 0$. This is carried out in Proposition \ref{prop:vmf} below. In the statement of this proposition, a tangent vector $U \in T_z\mathbb{R}^n$, at a point $z \in \mathbb{R}^n - \lbrace 0 \rbrace$, is written in the form
\begin{subequations}
\begin{equation} \label{eq:polarvector}
  U = u_\eta\,\bar{x} + \eta\,u
\end{equation}
where $z = \eta\,\bar{x}$, and where $u_\eta \in \mathbb{R}$ and $u \in T_{\bar{x}}S^{n-1}\,$. Here, $u_\eta$ and $u$ are unique, for a given $U$.  Precisely,
\begin{equation}
u_\eta \,=\, \langle U,\bar{x} \rangle \hspace{1cm} u \,=\, \frac{1}{\eta}\,\left[\, U \,-\, \langle U,\bar{x} \rangle\,\bar{x}\,\right]
\end{equation}
\end{subequations}
as follows since $\bar{x}$ and $u$ are orthogonal.
\begin{proposition} \label{prop:vmf}
the Rao-Fisher information metric $I$ of the von Mises-Fisher model is a well-defined Riemannian metric on the parameter space $\mathbb{R}^n$. On $\mathbb{R}^n - \lbrace 0 \rbrace$, it is a warped Riemannian metric of the form (\ref{eq:propexp1}), where
\begin{subequations} \label{eq:propvmf}
\begin{eqnarray} 
\label{eq:propvmf1}   \psi^{\prime\prime}(\eta) & = \frac{1}{n} + \frac{n-1}{n}\, \frac{I_{\nu+1}(\eta)}{\strut I_{\nu-1}(\eta)} - \frac{I^2_{\nu}(\eta)}{\strut I^2_{\nu-1}(\eta)}\\[0.12cm]
\label{eq:propvmf2}  \beta^2(\eta) & = \frac{\eta^2}{n} \,\left(\, 1 + \frac{I_{\nu+1}(\eta)}{\strut I_{\nu-1}(\eta)}\,\right) \hspace{1.1cm}
\end{eqnarray}
and it extends smoothly to the value,  
\begin{equation} \label{eq:propvmf3}
   I_0(U,U) \,=\, \frac{1}{n}\,\Vert U \Vert^2 
\end{equation}
at the point $z = 0$. Here, $\Vert\cdot\Vert$ denotes the Euclidean norm. 
\end{subequations} 
\end{proposition}
\vspace{0.1cm}

\textbf{Proof\,:} the Rao-Fisher information metric $I$ on $\mathbb{R}^n - \lbrace 0 \rbrace$ is given by Proposition \ref{prop:exponential}. This proposition applies because $M = S^{n-1}$ is an irreducible Riemannian symmetric space~\cite{helgason} (Table II, Page 354). Accordingly, for any point $z \in \mathbb{R}^n - \lbrace 0 \rbrace$, the metric $I_z$ is given by (\ref{eq:propexp1}). Formulae (\ref{eq:propvmf}) are proved as follows. \\[0.1cm]
\textit{Proof of (\ref{eq:propvmf1})\,:} this is carried out in Appendix \ref{app:B}, using the derivative and recurrence relations of modified Bessel functions~\cite{watson}.  \\[0.1cm] 
\textit{Proof of (\ref{eq:propvmf2})\,:} this follows from (\ref{eq:propexp2}) and (\ref{eq:vmf}). By (\ref{eq:vmf}), 
$$
\nabla_{\bar{x}}\,D(x,\bar{x}) \,=\, x - \langle x,\bar{x}\rangle\,\bar{x}
$$
which is just the orthogonal projection of $x$ onto the tangent space $T_{\bar{x}}S^{n-1\,}$. Replacing in (\ref{eq:propexp2}) gives
\begin{subequations} \label{eq:proofvmf}
\begin{equation} \label{eq:proofvmf1}
\beta^2(\eta) \,=\,\frac{\eta^2}{n-1}\,\mathbb{E}_z\left\Vert\nabla_{\bar{x}}\,D\right\Vert^2 \,=\,
\,\frac{\eta^2}{n-1}\,\mathbb{E}_z\left(\, 1 - \langle x,\bar{x}\rangle^2\,\right)
\end{equation}
Here, in the first equality, $n-1$ appears because $\mathrm{dim}\,S^{n-1} = n-1$. The second equality follows by Pythagoras' theorem, 
$$
\left\Vert\,x - \langle x,\bar{x}\rangle\,\bar{x}\,\right\Vert^2 \,=\,
\left\Vert\,x\,\right\Vert^2 - \left\Vert\,\langle x,\bar{x}\rangle\,\bar{x}\,\right\Vert^2 \,=\, 1 - \langle x,\bar{x}\rangle^2
$$ 
since $x$ and $\bar{x}$ belong to the unit sphere $S^{n-1\,}$. Formula (\ref{eq:propvmf2}) is derived from (\ref{eq:proofvmf1}) in Appendix \ref{app:B}, using the derivative and recurrence relations of modified Bessel functions~\cite{watson}. \\[0.1cm] 
\textit{Proof of (\ref{eq:propvmf3})\,:} for any point $z \in \mathbb{R}^n - \lbrace 0 \rbrace$, the metric $I_z$ is given by (\ref{eq:propexp1}). This reads
\begin{equation} \label{eq:proofvmf2}
   I_z(U,U) \,=\, \psi^{\prime\prime}(\eta)\,u^2_\eta \,+\, \beta^2(\eta)\, \Vert u\Vert^2
\end{equation}
Consider the limit of this expression at the point $z = 0$. In (\ref{eq:propvmf1}) and (\ref{eq:propvmf2}), this corresponds to the limit at $\eta = 0$. This can be evaluated using the power series development of modified Bessel functions~\cite{watson}. When replaced in (\ref{eq:propvmf1}) and (\ref{eq:propvmf2}), this gives the following developments,
$$
\begin{array}{llr}
\psi^{\prime\prime}(\eta)  = &\frac{1}{n} - &\hspace{-0.25cm}\frac{12}{\strut n^{\scriptscriptstyle 2}(n+2)}\,\left(\frac{\eta}{2}\right)^{\!2} + O\left(\eta^4\right) \\[0.14cm]
\nonumber \beta^2(\eta)  = & &\hspace{-0.4cm} \frac{\phantom{aa}4\phantom{aa}}{\mathstrut n}\,\left(\frac{\eta}{2}\right)^{\!2} + O\left(\eta^4\right)
\end{array}
$$
which immediately imply that
\begin{eqnarray}
\label{eq:proofvmf4}  \lim_{\eta\rightarrow 0}\,\, \psi^{\prime\prime}(\eta) & = \frac{1}{n} \\[0.1cm]
\label{eq:proofvmf5} \lim_{\eta\rightarrow 0}\,\, \beta^2(\eta) & = 0
\end{eqnarray}
Replacing (\ref{eq:proofvmf4}) and (\ref{eq:proofvmf5}) in (\ref{eq:proofvmf2}) gives,
\begin{equation} \label{eq:proofvmf9}
  \lim_{z\rightarrow 0}\,\, I_z(U,U) = \frac{1}{n}\,u^2_\eta
\end{equation}
Note that, from (\ref{eq:polarvector}),
$$
\Vert U \Vert^2 =  u^2_\eta \,+\, \eta^2\,\Vert u \Vert^2
$$
by Pythagoras' theorem, since $\bar{x}$ and $u$ are orthogonal. At the point $z = 0$, one has $\eta = 0$, so that
$$
 \Vert U \Vert^2_{z=0} =  u^2_\eta 
$$
This shows that (\ref{eq:proofvmf9}) is the same as
\begin{equation} \label{eq:proofvmf10}
  \lim_{z\rightarrow 0}\,\, I_z(U,U) = \frac{1}{n}\,\Vert U \Vert^2
\end{equation}
This limit does not depend on the path along which $z$ tends to $z = 0$. Therefore, $I_z$ extends smoothly to $I_{0\,}$, which is given by (\ref{eq:propvmf3}), at the point $z = 0$. This shows that $I$ is a well-defined Riemannian metric on the parameter space $\mathbb{R}^n$.
\end{subequations}
\hfill$\blacksquare$ 

\vfill
\pagebreak  

\subsection{The Riemannian Gaussian model} \label{subsec:gauss}
The Riemannian Gaussian model was recently introduced as a means of describing unimodal populations of covariance matrices~\cite{vemuri,said1,said2}. This model can be defined on any Riemannian symmetric space of non-positive sectional curvature. Let $M$ be such a symmetric space and denote $G$ its group of isometries. Then, the expressions appearing in (\ref{eq:exponential1}) are
\begin{equation} \label{eq:rgm}
   D(x,\bar{x}) = d^{\,2}(x,\bar{x}) \hspace{1cm} \eta(\sigma) = -\frac{1}{\strut 2\sigma^2}
\end{equation}
where $d(x,\bar{x})$ denotes the Riemannian distance in $M$, and condition (\ref{eq:exponential2}) is verified since each isometry $g \in G$ preserves this Riemannian distance. The function $\psi(\eta)$ is a strictly convex function of $\eta \in (-\infty\,,0)$, which can be expressed by means of a multiple integral~\cite{said2}, (see Proposition 1 in this reference). Precisely, $\psi(\eta)$ is the cumulant generating function of the squared Riemannian distance $d^{\,2}(x,\bar{x})$. 

Proposition \ref{prop:exponential} cannot be applied directly to the Riemannian Gaussian model (\ref{eq:rgm}). This is because, in most cases of interest, 
the Riemannian symmetric space $M$ is not irreducible. In such cases, before applying Proposition \ref{prop:exponential}, it is necessary to introduce the De Rham decomposition theorem~\cite{petersen}\cite{helgason}. \\[0.1cm]
\textbf{Remark 7\,:} assume the Riemannian symmetric space $M$ is moreover simply-connected. Then, the De Rham decomposition theorem implies that $M$ is a Riemannian product of irreducible Riemannian symmetric spaces~\cite{helgason} (Proposition 5.5, Page 310). Precisely, $M = M_1 \,\times\ldots\times\, M_r$ where each $M_q$ is an irreducible Riemannian symmetric space, and the Riemannian metric and distance of $M$ can be expressed as follows,
\begin{subequations} \label{eq:derham}
\begin{equation} \label{eq:derham1}
  Q(u,u) \,=\, \sum^r_{q=1}\, Q(u_{q\,},u_q)
\end{equation}
\begin{equation} \label{eq:derham2}
  d^{\,2}(x,y)\,=\, \sum^r_{q=1}\, d^{\,2}(x_{q\,},y_q)
\end{equation}
where $x,y \in M$ are written $x = (x_1\,,\ldots,\,x_r)$ and $y = (y_1\,,\ldots,\,y_r)$ with $x_{q\,},y_q \in M_{q\,}$, and where $u \in T_xM$ is written $u = u_1 \,+\ldots+\, u_r$ with $u_q \in T_{x_q}M_{q\,}$. Since $M$ has non-positive sectional curvature, each $M_q$ is either a Euclidean space, or a so-called space of non-compact type, having negative sectional curvature~\cite{helgason}. A concrete example of the De Rham decomposition is treated in Section \ref{sec:pn}, where $M = \mathcal{P}_n$ is the space of $n \times n$ real covariance matrices. \hfill$\blacksquare$ \\[0.1cm]
\end{subequations}
The following Proposition \ref{prop:rgm} gives the Rao-Fisher information metric of the Riemannian Gaussian model. Since this model is, in general, defined on a Riemannian symmetric space $M$ which is not irreducible, the Rao-Fisher information metric turns out to be a multiply-warped Riemannian metric, rather than a warped Riemannian metric. \\[0.1cm]
\textbf{Remark 8\,:}  in the notation of (\ref{eq:derham}), a multiply-warped Riemannian metric $I$ is a Riemannian metric defined on $\mathcal{M} \,= M\times (0\,,\infty)$, in the following way~\cite{multipl1}\cite{multipl2}
\begin{equation} \label{eq:mwarped}
  I_z(U,U) \,=\, \left(\,\!\alpha(\sigma) \, u_\sigma\,\!\right)^2 \,+\, \sum^r_{q=1}\beta^2_q(\sigma)\,Q_{\bar{x}}(u_{q\,},u_q) 
\end{equation}
for $z= (\bar{x},\sigma)$ and $U \in T_z\mathcal{M}$, where $U= u_\sigma\,\partial_\sigma\,+\,u$ with $u_\sigma \in \mathbb{R}$ and $u \in T_{\bar{x}}M$. Here, the functions $\alpha$ and $\beta_q$  are positive functions defined on $(0\,,\infty)$. Clearly, when $r = 1$, definition (\ref{eq:mwarped}) reduces to definition (\ref{eq:warpdef}) of a warped Riemannian metric. Therefore, warped Riemannian metrics are a special case of multiply-warped Riemannian metrics. \hfill$\blacksquare$ 
\begin{proposition} \label{prop:rgm}
 the Rao-Fisher information metric of the Riemannian Gaussian model is a multiply-warped Riemannian metric. In terms of $\bar{x} \in M$ and $\eta = -\left.1\middle/2\sigma^2\right.$, this metric has the following expression 
\begin{equation} \label{eq:proprgm1}
   I_z(U,U) \,=\,\psi^{\prime\prime}(\eta)\,u^2_\eta\,+\,\, \sum^r_{q=1}\,\left(4\eta^2\psi^\prime_q(\eta)\middle/\mathrm{dim}\,M_q\right)\,Q_{\bar{x}}(u_{q\,},u_q)
\end{equation}
where $U \,=\, u_\eta\,\partial_\eta\,+\,u$, and where $\psi_q(\eta)$ is the cumulant generating function of the squared Riemannian distance $d^{\,2}(x_{q\,},\bar{x}_q)$.
\end{proposition}
\vspace{0.1cm}
\textbf{Proof\,:} assume first that the Riemannian symmetric space $M$ is irreducible, so that Proposition \ref{prop:exponential} applies directly, and the Rao-Fisher information metric is given by (\ref{eq:propexp1}), 
\begin{subequations} \label{eq:irreducible}
\begin{equation} \label{eq:irreducible1}
   I_z(U,U) \,=\, \psi^{\prime\prime}(\eta)\,u^2_\eta\,+\, \beta^2(\eta)\,Q_{\bar{x}}(u,u) 
\end{equation}
To obtain $\beta^2(\eta)$, replace into (\ref{eq:propexp2}) the fact that
$$
\nabla_{\bar{x}}\,D(x,\bar{x}) \,=\, -2\,\mathrm{exp}^{-1}_{\bar{x}}(x) \hspace{1cm} Q\left(\nabla_{\bar{x}}\,D,\nabla_{\bar{x}}\,D\,\right) \,=\, 4\,d^{\,2}(x,\bar{x})
$$
where $\mathrm{exp}$ denotes the Riemannian exponential mapping, corresponding to the Riemannian metric $Q$ of $M$~\cite{chavel}, (see Page 407). It then follows from (\ref{eq:propexp2}) that,
\begin{equation} \label{eq:irreducible2}
\beta^2(\eta) \,=\, \left.4\eta^2\,\mathbb{E}_z\,d^{\,2}(x,\bar{x})\middle/\mathrm{dim}\,M\right. \,=\, \left.4\eta^2\,\psi^\prime(\eta)\middle/\mathrm{dim}\,M\right.
\end{equation}
where the second equality holds since $\psi(\eta)$ is the cumulant generating function of $d^{\,2}(x,\bar{x})$. From (\ref{eq:irreducible1}) and (\ref{eq:irreducible2}),
\begin{equation} \label{eq:irreducible3}
 I_z(U,U) \,=\, \psi^{\prime\prime}(\eta)\,u^2_\eta\,+\, \left(4\eta^2\psi^\prime(\eta)\middle/\mathrm{dim}\,M\right)\,Q_{\bar{x}}(u,u)
\end{equation}
\end{subequations}
which is the same as (\ref{eq:proprgm1}) with $r=1$. This proves the proposition in the special case where $M$ is irreducible. For the general case where $M$ is not irreducible, write $U \,=\, u_\eta\,\partial_\eta\,+\,u$, with $u = u_1 \,+\ldots+\, u_r$ as in Remark 7. It is possible to prove that,
\begin{subequations} \label{eq:reducible}
\begin{eqnarray} 
\label{eq:reducible2} q \neq p & \text{implies } & I_z(u_{p\,},u_q) =  0 \\[0.1cm]
\label{eq:reducible1} u = u_p &\!\!\text{ implies } & I_z(U,U)\,=\,\psi^{\prime\prime}(\eta)\,u^2_\eta\,+\, \left(4\eta^2\psi^\prime_p(\eta)\middle/\mathrm{dim}\,M_p\right)\,Q_{\bar{x}}(u_{p\,},u_p)
\end{eqnarray}
\end{subequations}
Then, since the Rao-Fisher information metric $I$ is bilinear and symmetric, (\ref{eq:proprgm1}) follows immediately, and the proposition is proved in the general case. 
\vfill
\pagebreak
\noindent \textit{Proof of identities (\ref{eq:reducible})\,:} this is carried out using the following properties (\ref{eq:gausspdf}). Note first that the probability density function of the Riemannian Gaussian model is given by (\ref{eq:exponential1}) and (\ref{eq:rgm}),
\begin{subequations} \label{eq:gausspdf}
\begin{equation} \label{eq:gausspdf1}
  p(x|\bar{x},\sigma) \,=\, \exp\left[\,\eta(\sigma)\,d^{\,2}(x,\bar{x}) \,-\, \psi(\eta(\sigma))\right]
\end{equation}
By substituting (\ref{eq:derham2}) in this expression, it is seen that
\begin{equation} \label{eq:gausspdf2}
  p(x|\bar{x},\sigma) \,=\, \prod^r_{q=1}\,\exp\left[\,\eta(\sigma)\,d^{\,2}(x_{q\,},\bar{x}_q) \,-\, \psi_q(\eta(\sigma))\right] \,=\,\prod^r_{q=1}\,p(x_q|\bar{x}_{q\,},\sigma)
\end{equation}
where $\psi_q(\eta)$ is the cumulant generating function of $d^{\,2}(x_{q\,},\bar{x}_q)$, as stated after (\ref{eq:proprgm1}). The last equality shows that $(x_q\,;q=1,\ldots,r)$ are independent, and that each $x_q$ has a Riemannian Gaussian density on the irreducible Riemannian symmetric space $M_{q\,}$, with parameters $z_q = (\bar{x}_{q\,},\sigma)$. Now, identities (\ref{eq:reducible}) can be obtained from definition (\ref{eq:rao}) of the Rao-Fisher information metric. To apply this definition, note from (\ref{eq:gausspdf2}), that the log-likelihood function $\ell(z)$ can be written,
\begin{equation} \label{eq:gausspdf3}
  \ell(z)(x) \,=\, \log p(x|z) \,=\,\sum^r_{q=1}\ell(z_q)(x_q) \hspace{1cm} \text{ where } \,\,\ell(z_q)(x_q) = \log p(x_q|z_q) \\[0.1cm]
\end{equation}
\end{subequations}
\textit{Proof of (\ref{eq:reducible2})\,:} recall the polarisation identity, from elementary linear algebra~\cite{lang}, (see Page 29),
$$
I_z(u_{p\,},u_q) \,=\, \frac{1}{4}\,I_z(u_p+u_{q\,},u_p+u_q) -  \frac{1}{4}\,I_z(u_p-u_{q\,},u_p-u_{q})
$$
By replacing (\ref{eq:rao}) into this identity, it can be seen that,
\begin{equation} \label{eq:ind1}
I_z(u_{p\,},u_q) \,=\, \mathbb{E}_z\left(\,\left(d\ell(z)\,u_p\right)\left(d\ell(z)\,u_q\right)\,\right)
\end{equation}
Using (\ref{eq:gausspdf3}), it is then possible to write
$$
I_z(u_{p\,},u_q) \,=\, \mathbb{E}_z\left(\,\left(d\ell(z)\,u_p\right)\left(d\ell(z)\,u_q\right)\,\right) \,=\,
\mathbb{E}_z\left(\,\left(d\ell(z_p)\,u_p\right)\left(d\ell(z_q)\,u_q\right)\,\right) \,=\, \mathbb{E}_{z_p}\left(d\ell(z_p)\,u_p\right)\,
\mathbb{E}_{z_q}\left(d\ell(z_q)\,u_q\right)
$$
Here, the first equality is just (\ref{eq:ind1}), while the second equality holds since $u_p \in T_{\bar{x}_p}M_p$ and $u_q \in T_{\bar{x}_q}M_q\,$, and the third equality holds since $x_p$ and $x_q$ are independent. Now, each one of the two expectations appearing on the right-hand side is equal to zero, since the expectation of the derivative of the log-likelihood must be zero~\cite{amari}, (see Page 28). This shows that (\ref{eq:reducible2}) holds. \hfill$\blacksquare$ \\[0.1cm]
\textit{Proof of (\ref{eq:reducible1})\,:} the condition $u = u_p$ implies $U = u_\eta\,\partial_\eta\,+\, u_{p\,}$. Replacing this in (\ref{eq:rao}), it follows using (\ref{eq:gausspdf3}),
\begin{subequations} \label{eq:prf}
\begin{equation} \label{eq:prf1}
I_z(U,U)\,=\, \mathbb{E}_z\left(\,\sum^r_{q=1} d\ell(z_q)\,U\,\right)^{\!2} = \mathbb{E}_z\left(\,\sum^r_{q=1} u_\eta\,\partial_\eta\ell(z_q) \,+\, d\ell(z_p)\,u_p\,\right)^{\!2}
\end{equation}
where the second equality holds since $u_p \in T_{\bar{x}_p}M_{p\,}$. Since the $x_q$ are independent, it is clear from (\ref{eq:gausspdf3}) that the $\ell(z_q)$ are independent. Accordingly, by expanding the right-hand side of (\ref{eq:prf1}),
\begin{equation} \label{eq:prf2}
  I_z(U,U)= \sum_{q\neq p} u^2_\eta\, \mathbb{E}_{z_q}\left(\partial_\eta\ell(z_q)\right)^2 \,+\, \mathbb{E}_{z_p}\left(\,u_\eta\,\partial_\eta \ell(z_p) \,+\, d\ell(z_p)\,u_p\,\right)^2 = \sum_{q\neq p} u^2_\eta\, \mathbb{E}_{z_q}\left(\partial_\eta\ell(z_q)\right)^2 + \mathbb{E}_{z_p}\left( \,d\ell(z_p)\,U\right)^2
\end{equation}
Applying (\ref{eq:amarari}) from the proof of Proposition \ref{prop:exponential} to each term in the sum over $q\neq p$, it follows that
\begin{equation} \label{eq:prf3}
I_z(U,U)= \sum_{q\neq p} \psi^{\prime\prime}_q(\eta)\, u^2_\eta \,+\, \mathbb{E}_{z_p}\left( \,d\ell(z_p)\,U\,\right)^2
\end{equation}
By (\ref{eq:rao}), the expectation appearing in the second term is given by the Rao-Fisher information metric of the Riemannian Gaussian model on the irreducible Riemannian symmetric space $M_{p\,}$. This can be replaced from (\ref{eq:irreducible3}), so that 
$$
\begin{array}{lll}
I_z(U,U) &= \sum_{q\neq p} \psi^{\prime\prime}_q(\eta)\, u^2_\eta \,+\, \psi^{\prime\prime}_p(\eta)\,u^2_\eta &+\, \left(4\eta^2\psi^\prime_p(\eta)\middle/\mathrm{dim}\,M_p\right)\,Q_{\bar{x}}(u_{p\,},u_p) \\[0.12cm]
                           &= \sum_{q} \psi^{\prime\prime}_q(\eta)\, u^2_\eta &+\, \left(4\eta^2\psi^\prime_p(\eta)\middle/\mathrm{dim}\,M_p\right)\,Q_{\bar{x}}(u_{p\,},u_p) 
\end{array}
$$
\end{subequations}
This immediately yields (\ref{eq:reducible1}), upon noting from (\ref{eq:gausspdf1}) and (\ref{eq:gausspdf2}) that $\psi(\eta) = \,{\scriptstyle \sum_q}\, \psi_q(\eta)$\,.
\hfill$\blacksquare$ \\[0.1cm]
Now, since identities (\ref{eq:reducible}) have been proved, (\ref{eq:proprgm1}) follows immediately, using the fact that the Rao-Fisher information metric $I$ is bilinear and symmetric. \hfill$\blacksquare$

\section{The generalised Mahalanobis distance} \label{sec:maha}
This section builds on Remark 3, made at the end of Section \ref{sec:warped}, in order to generalise the definition of the classical Mahalanobis distance, to the context of a location-scale model $\mathcal{P}$ defined on a Riemannian symmetric space $M$.

To begin, assume that, as in Theorem \ref{th:warplocdis}, the Riemannian symmetric space $M$ is irreducible and the location-scale model $\mathcal{P}$ verifies condition (\ref{eq:invariance}). Then, according to Theorem \ref{th:warplocdis}, the Rao-Fisher information metric $I$ of the model $\mathcal{P}$ is a warped Riemannian metric on the parameter space $\mathcal{M}$.

Recall from Remark 3 that this warped Riemannian metric $I$ induces an extrinsic Riemannian metric $Q^\sigma$ on $M$, for each $\sigma \in (0\,,\infty)$. The \textit{generalised Mahalanobis distance} is defined to be the Riemannian distance on $M$ which is induced by the extrinsic Riemannian metric $Q^\sigma$. The generalised Mahalanobis distance between $\bar{x}$ and $\bar{y}$ in $M$ is denoted $d(\bar{x},\bar{y}\,|\sigma)$. It is given by the following proposition.
\begin{proposition} \label{prop:maha1}
 the generalised Mahalanobis distance $d(\bar{x},\bar{y}\,|\sigma)$ between $\bar{x}$ and $\bar{y}$ in $M$ is given by
\begin{equation} \label{eq:mahalanobis}
  d(\bar{x},\bar{y}\,|\sigma) \,=\, \beta(\sigma)\,d(\bar{x},\bar{y})
\end{equation}
where the function $\beta(\sigma)$ is given by (\ref{eq:warplocdis}), and where $d(\bar{x},\bar{y})$ denotes the Riemannian distance in $M$. 
\end{proposition}
\vfill
\pagebreak
\textbf{Remark 9\,:} the generalised Mahalanobis distance (\ref{eq:mahalanobis}) reduces to the classical Mahalanobis distance, when $\mathcal{P}$ is the isotropic normal model on $M = \mathbb{R}^d$. In this case, the Rao-Fisher metric $I$ is given by (\ref{eq:poincare1}) in the introduction, so that $\beta(\sigma) = 1/\sigma$. Replacing this in (\ref{eq:mahalanobis}) yields
\begin{equation} \label{eq:mahalanobisr}
  d(\bar{x},\bar{y}\,|\sigma) \,=\, \frac{1}{\sigma}\,\Vert\bar{x} - \bar{y}\Vert 
\end{equation}
where $\Vert\bar{x} - \bar{y}\Vert$ is the Euclidean distance in $M = \mathbb{R}^d$. Now, (\ref{eq:mahalanobisr}) is the classical Mahalanobis distance~\cite{maclahlan}. 
\hfill$\blacksquare$ \\[0.1cm]
\textbf{Proof of Proposition \ref{prop:maha1}\,:} the extrinsic metric $Q^\sigma$ is induced by the warped Riemannian metric $I$, which is given by (\ref{eq:warplocdis}). Therefore, it follows from (\ref{eq:inducedmetric}) that
\begin{subequations}
\begin{equation} \label{eq:inducedmetric1}
 Q^\sigma_{\bar{x}}(u,u) \,=\, \beta^2(\sigma)\,Q_{\bar{x}}(u,u) 
\end{equation}
where the function $\beta(\sigma)$ is given by (\ref{eq:warplocdis}). To find the generalised Mahalanobis distance between $\bar{x}$ and $\bar{y}$ in $M$, let $c(t)$ be a curve in $M$ with $c(0) = \bar{x}$ and $c(1) = \bar{y}$. Denote by $L(c|\sigma)$ and by $L(c)$ the length of this curve, as measured by the Riemannian metrics $Q^\sigma$ and $Q$, respectively. Then, 
\begin{equation} \label{eq:inducedmetric2}
L(c|\sigma) \,=\, \int^1_0\, \left(Q^\sigma(\dot{c},\dot{c})\,\right)^{1/2}\,dt \,=\, \beta(\sigma)\, \int^1_0\, \left(Q(\dot{c},\dot{c})\,\right)^{1/2}\,dt \,=\, \beta(\sigma)\,L(c)
\end{equation}
where the second equality follows from (\ref{eq:inducedmetric1}). To obtain (\ref{eq:mahalanobis}), it is enough to write
\begin{equation} \label{eq:inducedmetric3}
d(\bar{x},\bar{y}\,|\sigma) \,=\, \inf_c\, L(c|\sigma) \,=\, \beta(\sigma)\, \inf_c\,L(c) \,=\, \beta(\sigma)\,d(\bar{x},\bar{y})
\end{equation}
where the infimum is over all curves $c(t)$ as above, and the second equality follows from (\ref{eq:inducedmetric2}). \hfill$\blacksquare$ \\[0.1cm]
\end{subequations}
Expression (\ref{eq:mahalanobis}) of the generalised Mahalanobis distance is valid only under the assumption that the Riemannian symmetric space $M$ is irreducible. This assumption does not hold, when the model $\mathcal{P}$ is the Riemannian Gaussian model studied in \ref{subsec:gauss}. For this model, an alternative expression of the generalised Mahalanobis distance is given in Proposition \ref{prop:maha2} below.

As in \ref{subsec:gauss}, let $\mathcal{P}$ be the Riemannian Gaussian model on a Riemannian symmetric space $M$, where $M$ is simply-connected and has non-positive sectional curvature. Proposition \ref{prop:rgm} states that the Rao-Fisher information metric $I$ of the model $\mathcal{P}$ is a multiply-warped Riemannian metric on the parameter space $\mathcal{M}$. For each $\sigma \in (0\,,\infty)$, this multiply-warped Riemannian metric $I$ induces an extrinsic Riemannian metric $Q^\sigma$ on $M$. Precisely, $Q^\sigma$ can be obtained from (\ref{eq:proprgm1}),
\begin{equation} \label{eq:rgmmaha1}
 Q^\sigma_{\bar{x}}(u,u) \,= \, \sum^r_{q=1}\,\beta^2_q(\sigma)\,Q_{\bar{x}}(u_{q\,},u_q) \hspace{1cm}
\beta^2_q(\sigma) \,=\, \left.4\eta^2\psi^\prime_q(\eta)\middle/\mathrm{dim}\,M_q\right.
\end{equation}
The generalised Mahalanobis distance $d(\bar{x},\bar{y}\,|\sigma)$ is the Riemannian distance between $\bar{x}$ and $\bar{y}$ in $M$, induced by this extrinsic Riemannian metric $Q^\sigma$. 
\begin{proposition} \label{prop:maha2}
when $\mathcal{P}$ is the Riemannian Gaussian model, the generalised Mahalanobis distance $d(\bar{x},\bar{y}\,|\sigma)$ between $\bar{x}$ and $\bar{y}$ in $M$ is given by
\begin{equation} \label{eq:mahalanobisrgm} 
  d^2(\bar{x},\bar{y}\,|\sigma) \,=\, \sum^r_{q=1}\,\beta^2_q(\sigma)\,d^{\,2}(\bar{x}_{q\,},\bar{y}_q)
\end{equation}
where the notation is that of (\ref{eq:derham2}).
\end{proposition}
\textbf{Proof\,:} the proof hinges on the fact that the extrinsic Riemannian metric $Q^\sigma$ of (\ref{eq:rgmmaha1}) is an invariant Riemannian metric on $M$. In other words, if $G$ is the group of isometries of $M$, then
\begin{subequations} 
\begin{equation} \label{eq:rgmmaha2}
Q^\sigma_{g\cdot\bar{x}}(dg_{\bar{x}}\,u,dg_{\bar{x}}\,u) = Q^\sigma_{\bar{x}}(u,u) \hspace{1cm} \text{ for all } g \in G
\end{equation}
where $dg_{\bar{x}}$ is the derivative of the isometry $g$ at the point $\bar{x}$. The proof of (\ref{eq:rgmmaha2}) is not detailed here. It follows since the Riemannian metric $Q$ is also an invariant Riemannian metric on $M$, so that $Q$ also verifies (\ref{eq:rgmmaha2}), and since $Q^\sigma$ is related to $Q$ by (\ref{eq:rgmmaha1}). A general result in~\cite{helgason} (Corollary 4.3, Page 182), states that all invariant Riemannian metrics on $M$ have the same geodesics. In particular, the metrics $Q^\sigma$ and $Q$ have the same geodesics, and therefore the same Riemannian exponential mapping $\mathrm{exp}$.  To find the generalised Mahalanobis distance between $\bar{x}$ and $\bar{y}$ in $M$, let $u = \mathrm{exp}^{-1}_{\bar{x}}(\bar{y})$, and note that
\begin{equation} \label{eq:invg1}
d^2(\bar{x},\bar{y}\,|\sigma) \,=\, Q^\sigma_{\bar{x}}(u,u) \,=\,  
\sum^r_{q=1}\,\beta^2_q(\sigma)\,Q_{\bar{x}}(u_{q\,},u_{q})
\end{equation}
where the second equality follows from (\ref{eq:rgmmaha1}). Now, to prove (\ref{eq:mahalanobisrgm}) it is enough to prove that
\begin{equation} \label{eq:invg2}
Q_{\bar{x}}(u_{q\,},u_{q}) \,=\,d^{\,2}(\bar{x}_{q\,},\bar{y}_q)
\end{equation}
Indeed, (\ref{eq:mahalanobisrgm}) is then obtained by replacing (\ref{eq:invg2}) into (\ref{eq:invg1}). The proof of (\ref{eq:invg2}) follows by writing, as in (\ref{eq:invg1}),
\begin{equation} \label{eq:invg3}
d^2(\bar{x},\bar{y}) \,=\, Q_{\bar{x}}(u,u) \,=\,  
\sum^r_{q=1}\,Q_{\bar{x}}(u_{q\,},u_{q}) \,=\, \sum^r_{q=1}\, d^{\,2}(\bar{x}_{q\,},\bar{y}_q)
\end{equation}
where the second equality follows from (\ref{eq:derham1}), and the third equality follows from (\ref{eq:derham2}). Since (\ref{eq:invg3}) is an identity which holds for arbitrary $\bar{x} = (\bar{x}_1\,,\ldots,\,\bar{x}_r)$ and $\bar{y} = (\bar{y}_1\,,\ldots,\,\bar{y}_r)$, it follows that (\ref{eq:invg2}) must hold true, as required. 
\end{subequations}
\hfill $\blacksquare$ \\[0.1cm]
\textbf{Remark 10\,:} the generalised Mahalanobis distance, whether given by (\ref{eq:mahalanobis}) or by (\ref{eq:mahalanobisrgm}), is an invariant Riemannian distance on $M$,
\begin{equation} \label{eq:invmaha}
d(g\cdot\bar{x},g\cdot\bar{y}\,|\sigma) \,= \, d(\bar{x},\bar{y}\,|\sigma) \hspace{1cm} \text{ for all } g \in G
\end{equation}
This follows since the Riemannian distance in $M$ is invariant\,: $d(g\cdot\bar{x},g\cdot\bar{y}) \,= \, d(\bar{x},\bar{y})$. \hfill$\blacksquare$
\section{A concrete example for the Riemannian Gaussian model} \label{sec:pn}
The aim of this section is to illustrate the geometric concepts involved in Propositions \ref{prop:rgm} and \ref{prop:maha2}, by applying these concepts to the concrete example of the Riemannian Gaussian model defined on $M = \mathcal{P}_{n\,}$, the space of $n \times n$ real covariance matrices.

The space $\mathcal{P}_{n\,}$ is a Riemannian symmetric space, which is simply-connected and has non-positive sectional curvature~\cite{helgason}\cite{terras2}. It is usually equipped with its affine-invariant Riemannian metric~\cite{terras2}\cite{atkinson}, 
\begin{subequations} \label{eq:affine}
\begin{equation} \label{eq:affine1}
  Q_{\bar{x}}(u,u) \,=\, \mathrm{tr}\left[\bar{x}^{-1}u\right]^2 \hspace{1cm} \bar{x} \in \mathcal{P}_n \,,\, u \in T_{\bar{x}}\mathcal{P}_n
\end{equation}
This metric is invariant under the action of the group of isometries $G = GL(n,\mathbb{R})$ on $\mathcal{P}_{n\,}$, which is given by affine transformations,
\begin{equation} \label{eq:affine2}
  g\cdot \bar{x} = g\,\bar{x}\,g^t
\end{equation}
where $^t$ denotes the transpose. Moreover, this metric induces a Riemannian distance on $\mathcal{P}_{n\,}$, which is given by, 
\begin{equation} \label{eq:affine3}
  d^2(\bar{x},\bar{y}) \,=\, \mathrm{tr} \left[\log\left(\bar{x}^{-1/2}\,\bar{y}\,\bar{x}^{-1/2}\right)\right]^2
\end{equation}
This distance is also invariant under the action of the group $GL(n,\mathbb{R})$ on $\mathcal{P}_{n\,}$. In other words\,: $d(g\cdot\bar{x},g\cdot\bar{y}) \,= \, d(\bar{x},\bar{y})$. 
\end{subequations}

\begin{subequations} \label{eq:rgmpn}
The Riemannian Gaussian model on $\mathcal{P}_n\,$ is given by the probability density function~\cite{vemuri}\cite{said1} 
\begin{equation} \label{eq:rgmpn1}
p(x|\bar{x},\sigma) \,=\, Z^{-1}(\sigma)\,\exp\left[\,-\frac{d^2(x,\bar{x})}{2\sigma^2}\,\right]
\end{equation}
which is a probability density function with respect to the invariant volume element associated to the Riemannian metric (\ref{eq:affine1}). The normalising factor $Z(\sigma)$ can be expressed as a multiple integral~\cite{said1}, (see Proposition 4 in this reference), 
\begin{equation} \label{eq:rgmpn2}
  Z(\sigma) \,=\, \mathrm{C}_n\,\int_{\mathbb{R}^n}\, e^{-\left.\Vert r \Vert^2\middle/ 2\sigma^2\right.}\, \prod_{i <j}\,\sinh\left(|r_i-r_j|\middle/2\right)\,dr_1\,\ldots\,dr_n
\end{equation}
where $\mathrm{C}_n$ is a numerical constant which only depends on $n$, and the integration variable is denoted $r = (r_1\,,\ldots,\,r_n) \in \mathbb{R}^n$. If $\eta(\sigma) = \left.-1\middle/2\sigma^2\right.$, then $\psi(\eta) \,=\, \log\,Z(\sigma)$ is a strictly convex function of $\eta \in (-\infty\,,0)$.
\end{subequations}

\begin{subequations} \label{eq:derhampn}
With (\ref{eq:affine}) and (\ref{eq:rgmpn}) in mind, consider the application of Proposition \ref{prop:rgm} to the Riemannian Gaussian model on $\mathcal{P}_{n\,}$. This will lead to the expression of the Rao-Fisher information metric $I$ of this model. 

\textbf{De Rham decomposition of $\mathcal{P}_{n}$\,:} recall first that Proposition \ref{prop:rgm} uses the De Rham decomposition, introduced in Remark 7. For the Riemannian symmetric space $\mathcal{P}_{n\,}$, the De Rham decomposition states that $\mathcal{P}_{n}$ is a Riemannian product of irreducible Riemannian symmetric spaces $\mathcal{P}_n = \mathbb{R} \times S\mathcal{P}_n\,$, where $S\mathcal{P}_n$ is the set of $\bar{s} \in \mathcal{P}_n$ such that $\det(\bar{s}) = 1$. The identification of $\mathcal{P}_n$ with $\mathbb{R} \times S\mathcal{P}_n$ is obtained by identifying each $\bar{x} \in \mathcal{P}_n$ with a couple $(\bar{\tau},\bar{s})$, where $\bar{\tau} \in \mathbb{R}$ and $\bar{s} \in S\mathcal{P}_n$ are given by
\begin{equation} \label{eq:derhampn1}
   \bar{\tau} = \log\det(\bar{x}) \hspace{1cm} \bar{s} = e^{-\left.\bar{\tau}\middle/n\right.}\,\bar{x}
\end{equation}
Note that the spaces $\mathbb{R}$ and $S\mathcal{P}_n$ are indeed irreducible Riemannian symmetric spaces. This is clear for $\mathbb{R}$, which is one-dimensional and cannot be decomposed into a product of lower-dimensional spaces. The fact that $S\mathcal{P}_n$ is irreducible can be found in~\cite{helgason}, (Table II, Page 354). It will be convenient to write $\bar{x} = (\bar{x}_1,\bar{x}_2)$ where $\bar{x}_1 = \bar{\tau}$ and $\bar{x}_2 = \bar{s}$. If $u \in T_{\bar{x}}\mathcal{P}_n$, then $u = u_1 \,+\, u_2\,$, 
\begin{equation} \label{eq:derhampn2}
  u_1 = \frac{1}{n}\,\mathrm{tr}(\bar{x}^{-1}u)\,\bar{x} \hspace{1cm} u_2 = u - \frac{1}{n}\,\mathrm{tr}(\bar{x}^{-1}u)\,\bar{x}
\end{equation}
Here, $u_1 \in T_{\bar{x}_1}\mathbb{R}$, where $T_{\bar{x}_1}\mathbb{R} \subset T_{\bar{x}}\mathcal{P}_n$ is the one-dimensional subspace consisting of symmetric matrices $v$ of the form $v = t\,\bar{x}$ with $t$ any real number. On the other hand, $u_2\in T_{\bar{x}_2}S\mathcal{P}_n$, where $T_{\bar{x}_2}S\mathcal{P}_n \subset T_{\bar{x}}\mathcal{P}_n$ is the subspace consisting of symmetric matrices $v$ which satisfy $\mathrm{tr}(\bar{x}^{-1}v) = 0$. 
Using (\ref{eq:derhampn1}) and (\ref{eq:derhampn2}), (\ref{eq:derham1}) and (\ref{eq:derham2}) of Remark 7 can be written down, 
\begin{equation} \label{eq:derhampn3}
 Q_{\bar{x}}(u,u) \,=\, Q_{\bar{x}}(u_{1\,},u_1) + Q_{\bar{x}}(u_{2\,},u_2)
\end{equation}
\begin{equation} \label{eq:derhampn4}
  d^2(\bar{x},\bar{y}) \,=\, \frac{1}{n}\,\left|\bar{x}_1 - \bar{y}_1\right|^2 + d^2(\bar{x}_{2\,},\bar{y}_2)
\end{equation}
where $Q_{\bar{x}}$ is the affine-invariant metric (\ref{eq:affine1}) and $d(\bar{x},\bar{y})$ or $d(\bar{x}_{2\,},\bar{y}_2)$ is the Riemannian distance (\ref{eq:affine3}). The proof of formulae (\ref{eq:derhampn3}) and (\ref{eq:derhampn4}) is a direct calculation, and is not detailed here.
\end{subequations}

\begin{subequations} \label{eq:raopn}
\textbf{The Rao-Fisher metric $I$\,:} according to (\ref{eq:proprgm1}) of Proposition \ref{prop:rgm}, the Rao-Fisher information metric $I$ of the Riemannian Gaussian model on $\mathcal{P}_n$ is given by,
\begin{equation} \label{eq:raopn1}
  I_z(U,U) = \psi^{\prime\prime}(\eta)\,u^2_\eta\,+\, \left(4\eta^2\psi^\prime_1(\eta)\middle/\mathrm{dim}\,\mathbb{R}\right)\,Q_{\bar{x}}(u_{1\,},u_1)\,+\,
 \left(4\eta^2\psi^\prime_2(\eta)\middle/\mathrm{dim}\,S\mathcal{P}_n\right)\,Q_{\bar{x}}(u_{2\,},u_2)
\end{equation}
for $z = (\bar{x},\sigma)$ in the parameter space $\mathcal{M} = \mathcal{P}_n \times (0\,,\infty)$, and for $U = u_\eta\,\partial_\eta\,+\, u$ where $u = u_1 + u_2$ is given by (\ref{eq:derhampn2}). Indeed, (\ref{eq:raopn1}) results from (\ref{eq:proprgm1}), by putting $r = 2$, as well as $M_1 = \mathbb{R}$ and $M_2 = S\mathcal{P}_n$. The functions appearing in (\ref{eq:raopn1}) are $\psi(\eta) = \log\,Z(\sigma)$ with $Z(\sigma)$ given by (\ref{eq:rgmpn2}), and, as shown in Remark 11 below,
\begin{equation} \label{eq:raopn2}
  \psi_1(\eta) = \frac{1}{2}\,\log(2\pi n) - \frac{1}{2}\,\log\left(-2\eta\right) \hspace{1cm} \psi_2(\eta) = \psi(\eta) - \psi_1(\eta)
\end{equation}
Moreover, $\mathrm{dim}\,\mathbb{R} = 1$ and $\mathrm{dim}\,S\mathcal{P}_n = \mathrm{dim}\,\mathcal{P}_n - 1 = n(n+1)/2 - 1$. Replacing into (\ref{eq:raopn1}) gives,
\begin{equation} \label{eq:raopn3}
  I_z(U,U) = \psi^{\prime\prime}(\eta)\,u^2_\eta\,-\, 2\eta\,Q_{\bar{x}}(u_{1\,},u_1)\,+\,
 \frac{8\eta^2\psi^\prime_2(\eta)}{\mathstrut n^2+n-2}\,Q_{\bar{x}}(u_{2\,},u_2)
\end{equation}
This expression of the Rao-Fisher information metric of the Riemannian Gaussian model on $\mathcal{P}_n$ can be computed directly from (\ref{eq:affine1}), (\ref{eq:derhampn2}) and (\ref{eq:raopn2}), once the function $\psi(\eta)$ is known. This function $\psi(\eta)$ has been tabulated for values of $n$ up to $n = 50$, using a Monte Carlo method which was developed specifically for the evaluation of (\ref{eq:rgmpn2})~\cite{paolo}.
\vfill
\pagebreak
\noindent \textbf{Remark 11\,:} assume $x$ follows the Riemannian Gaussian probability density (\ref{eq:rgmpn1}) on $\mathcal{P}_n\,$. If $x = (x_{1\,},x_2)$ where $x_1 \in \mathbb{R}$ and $x_2 \in S\mathcal{P}_n\,$, then the densities of $x_1$ and $x_2$ can be found by replacing (\ref{eq:derhampn4}) into (\ref{eq:rgmpn1}). Precisely, this gives
$$
p(x|\bar{x},\sigma) \,\propto\, \exp\left[\, - \frac{\left|x_1-\bar{x}_1\right|^2}{2n\sigma^2}\,\right]\,\times\,
\exp\left[\,-\frac{d^2(x_{2\,},\bar{x}_2)}{2\sigma^2}\,\right]
$$
It follows from this decomposition that $x_1$ and $x_2$ are independent, and that $x_1$ follows a univariate normal distribution of mean $\bar{x}_1$ and of variance $n\sigma^2$. In particular, the moment generating function $\psi_1(\eta)$ of the squared distance $\left|x_1-\bar{x}_1\right|^2$ has the expression stated in (\ref{eq:raopn2}). \hfill$\blacksquare$ \\[0.1cm]
\indent \textbf{The generalised Mahalanobis distance on $\mathcal{P}_n$\,:} applying Proposition \ref{prop:maha2} will yield the expression of the generalised Mahalanobis distance on $\mathcal{P}_n$. The Rao-Fisher information metric $I$ as given by (\ref{eq:raopn3}) induces an extrinsic Riemannian metric $Q^\sigma$ on $\mathcal{P}_n\,$, for each $\sigma \in (0\,,\infty)$,
\end{subequations}
\begin{subequations} \label{eq:mahapn}
\begin{equation} \label{eq:mahapn1}
  Q^\sigma_{\bar{x}}(u,u) \,=\, -\, 2\eta\,Q_{\bar{x}}(u_{1\,},u_1)\,+\,
 \frac{8\eta^2\psi^\prime_2(\eta)}{\mathstrut n^2+n-2}\,Q_{\bar{x}}(u_{2\,},u_2) \hspace{1cm} \eta = -\frac{1}{\mathstrut 2\sigma^2}
\end{equation}
The generalised Mahalanobis distance on $\mathcal{P}_n$ is the Riemannian distance induced on $\mathcal{P}_n$ by the extrinsic Riemannian metric $Q^\sigma$. If the generalised Mahalanobis distance between $\bar{x}$ and $\bar{y}$ in $\mathcal{P}_n$ is denoted $d(\bar{x},\bar{y}\,|\sigma)$, then (\ref{eq:mahalanobisrgm}) of  Proposition \ref{prop:maha2}, along with (\ref{eq:mahapn1}), imply
\begin{equation} \label{eq:mahapn2} 
  d^2(\bar{x},\bar{y}\,|\sigma) = \frac{\left|\bar{x}_1 - \bar{y}_1\right|^2}{\mathstrut n\sigma^2} \,+\, \frac{4\psi^\prime_2\left(-\frac{1}{\mathstrut 2\sigma^2}\right)}{\mathstrut (n^2+n-2)\sigma^4}\,d^2(\bar{x}_{2\,},\bar{y}_2)
\end{equation}
This distance can be computed directly from (\ref{eq:affine3}), (\ref{eq:derhampn1}) and (\ref{eq:raopn2}), once the function $\psi(\eta)$ has been tabulated using the Monte Carlo method of~\cite{paolo}, or computed in any other way. 
\end{subequations}

\begin{subequations} \label{eq:affinv}
\textbf{Affine invariance of the generalised Mahalanobis distance\,:} the affine-invariant Riemannian metric $Q$ of (\ref{eq:affine1}) is well-known to the information science community, having been introduced in~\cite{pennec2}. Besides the metric $Q$, a whole new family of affine-invariant Riemannian metrics $Q^\sigma$ is provided by (\ref{eq:mahapn1}). Indeed, to say that $Q$ is affine-invariant means that it is invariant under affine transformations (\ref{eq:affine2}). In other words
\begin{equation} \label{eq:affinv1}
  Q_{g\cdot \bar{x}}(dg_{\bar{x}}\,u,dg_{\bar{x}}\,u) \,=\, Q_{\bar{x}}(u,u) \hspace{1cm} \text{ for all } g \in GL(n,\mathbb{R})
\end{equation}
where $dg_{\bar{x}}$ denotes the derivative of the affine transformation (\ref{eq:affine2}) at the point $\bar{x} \in \mathcal{P}_n\,$. On the other hand, it is shown in Remark 12 below that each one of the metrics $Q^\sigma$ also verifies (\ref{eq:affinv1}), so that
\begin{equation} \label{eq:affinv2}
  Q^\sigma_{g\cdot \bar{x}}(dg_{\bar{x}}\,u,dg_{\bar{x}}\,u) \,=\, Q^\sigma_{\bar{x}}(u,u) \hspace{1cm} \text{ for all } g \in GL(n,\mathbb{R})
\end{equation}
This means that each one of the metrics $Q^\sigma$ is an affine-invariant Riemannian metric, as claimed. Furthermore, the fact that the metric $Q^\sigma$ is invariant under affine transformations implies that the generalised Mahalanobis distance (\ref{eq:mahapn2}) is also invariant under these transformations,
\begin{equation} \label{eq:affinv3}
  d(g\cdot\bar{x},g\cdot\bar{y}\,|\sigma) \,=\, d(\bar{x},\bar{y}\,|\sigma) \hspace{1cm} \text{ for all } g \in GL(n,\mathbb{R})
\end{equation}
\end{subequations}
This is because the generalised Mahalanobis distance (\ref{eq:mahapn2}) is the Riemannian distance induced on $\mathcal{P}_n$ by $Q^\sigma$. \\[0.1cm]
\noindent \textbf{Remark 12\,:} to prove that (\ref{eq:affinv2}) holds for each one of the metrics $Q^\sigma$, write (\ref{eq:mahapn1}) in the form
\begin{subequations}
\begin{equation} \label{eq:rk121}
  Q^\sigma_{\bar{x}}(u,u) \,=\, \beta^2_1(\sigma)\,Q_{\bar{x}}(u_{1\,},u_1) \,+\, \beta^2_2(\sigma)\,Q_{\bar{x}}(u_{2\,},u_2)
\end{equation}
Let $g \in GL(n,\mathbb{R})$ and let $\bar{y} = g\cdot \bar{x}$ and $v = dg_{\bar{x}}\,u$, so that $v \in T_{\bar{y}}\mathcal{P}_n\,$. Using (\ref{eq:rk121}), the left-hand side of (\ref{eq:affinv2}) takes on the form
\begin{equation} \label{eq:rk122}
Q^\sigma_{\bar{y}}(v,v) \,=\, \beta^2_1(\sigma)\,Q_{\bar{y}}(v_{1\,},v_1) \,+\, \beta^2_2(\sigma)\,Q_{\bar{y}}(v_{2\,},v_2)
\end{equation}
Note from (\ref{eq:affine2}) that
\begin{equation} \label{eq:rk123}
\bar{y} = g\cdot \bar{x}\,=\,g\,\bar{x}\,g^t \hspace{1cm} v = dg_{\bar{x}}\,u \,=\, g\,u\,g^t
\end{equation}
Therefore, it follows from (\ref{eq:derhampn2}) that $v_1$ is given by
$$
v_1 = \frac{1}{n}\,\mathrm{tr}\left[\bar{y}^{-1}v\right]\,\bar{y} \,=\, \frac{1}{n}\,\mathrm{tr}\left[\bar{x}^{-1}u\right]\,\bar{y} 
$$
where the second equality follows by a direct calculation. From (\ref{eq:rk123}), it is now seen that
\begin{equation} \label{eq:boo}
v_1 = \frac{1}{n}\,\mathrm{tr}\left[\bar{x}^{-1}u\right]\,g\,\bar{x}\,g^t = g\,u_1\,g^t = dg_{\bar{x}}\,u_1 
\end{equation}
where the second equality follows from (\ref{eq:derhampn2}). Since $v = v_1 + v_2\,$, this implies that $v_2$ is given by
\begin{equation} \label{eq:baa}
v_2 \,=\, v - v_1 \,=\, g\,u\,g^t - g\,u_1\,g^t \,= \, g\,(u - u_1)\,g^t \,=\, g\,u_2\,g^t \,=\, dg_{\bar{x}}\,u_2
\end{equation}
where the second equality follows from (\ref{eq:rk123}) and (\ref{eq:boo}). Using (\ref{eq:affinv1}), it follows from (\ref{eq:boo}) and (\ref{eq:baa}) that
\begin{equation} \label{eq:rk124}
Q_{\bar{y}}(v_{1\,},v_1) = Q_{\bar{x}}(u_{1\,},u_1)  \hspace{1cm} Q_{\bar{y}}(v_{2\,},v_2) = Q_{\bar{x}}(u_{2\,},u_2)
\end{equation}
Finally, replacing (\ref{eq:rk124}) into (\ref{eq:rk122}), it is found that 
$$
Q^\sigma_{\bar{y}}(v,v) \,=\, \beta^2_1(\sigma)\,Q_{\bar{x}}(u_{1\,},u_1) \,+\, \beta^2_2(\sigma)\,Q_{\bar{x}}(u_{2\,},u_2) \,=\, Q^\sigma_{\bar{x}}(u,u)
$$
where the second equality follows from (\ref{eq:rk121}). Recalling the definitions of $\bar{y}$ and $v$, it is clear that the proof of (\ref{eq:affinv2}) is now complete. \hfill$\blacksquare$
\end{subequations}
\vfill
\pagebreak


\section{The solution of the geodesic equation} \label{sec:geodesic}
The present section provides the solution of the geodesic equation of a multiply-warped Riemannian metric. The main result is the following Proposition \ref{prop:geodesicmultwarp}. This proposition shows that the solution of the geodesic equation of a multiply-warped Riemannian metric, for given initial conditions, reduces to the solution of a one-dimensional second-order differential equation. As stated in Remark 8, warped Riemannian metrics are a special case of multiply-warped Riemannian metrics. Therefore, Proposition \ref{prop:geodesicmultwarp} also applies to the solution of the geodesic equation of a warped Riemannian metric. This special case of warped Riemannian metrics was treated separately in~\cite{oneil}.

Let $I$ be a multiply-warped Riemannian metric defined on $\mathcal{M} \,=\, M\times (0\,,\infty)$, in the notation of (\ref{eq:mwarped}),
\begin{equation} \label{eq:mwarpedbis}
  I_z(U,U) \,=\, \left(\,\!\alpha(\sigma) \, u_\sigma\,\!\right)^2 \,+\, \sum^r_{q=1}\beta^2_q(\sigma)\,Q_{\bar{x}}(u_{q\,},u_q) 
\end{equation}
for $z = (\bar{x},\sigma)$ and $U \in T_z\mathcal{M}$, with $u = u_\sigma\,\partial_\sigma\,+\,u$ and $u = u_1\,+\ldots+\,u_r\,$. As in (\ref{eq:warpdefr}) of Section \ref{sec:warped}, introduce the vertical distance coordinate $r$, which is defined by $dr/d\sigma = \alpha(\sigma)$.
\vspace{0.1cm} 
\begin{proposition} \label{prop:geodesicmultwarp}
let $\gamma(t)$ be a geodesic of the multiply-warped Riemannian metric $I$, with initial conditions $\gamma(0) = z$ and $\dot{\gamma}(0) = U$, and let $\gamma(t) = (\bar{x}(t),\sigma(t))$ and $r(t) = r(\sigma(t))$. Then, $r(t)$ verifies the second-order differential equation
\begin{subequations} \label{eq:geodesic}
\begin{equation} \label{eq:geodesic1}
  \ddot{r} \,=\, -\frac{1}{2}\frac{d}{dr}V(r) \hspace{1cm} V(r) \,=\, \sum^r_{q=1}\,\frac{\beta^2_q(r(0))}{\beta^2_q(r)}\,I_z(u_{q\,},u_q)
\end{equation}
and $\bar{x}(t)$ is given by
\begin{equation} \label{eq:geodesic2}
  \bar{x}(t) \,=\, \exp_{\bar{x}}\left[\,\sum^r_{q=1}\left(\,\int^{\scriptscriptstyle t}_{\scriptscriptstyle 0}\frac{\beta^2_q(r(0))}{\strut\beta^2_q(r(s))}ds\,\right) \,u_q\,\right]
\end{equation}
where $\exp$ denotes the Riemannian exponential mapping of the metric $Q$ on $M$.
\end{subequations}
\end{proposition}
\textbf{Proof\,:} the proof is given in Appendix \ref{app:C}. It is a generalisation of the proof dealing with the special case of warped Riemannian metrics, which can be found in~\cite{oneil} (Proposition 38, Page 208). \hfill$\blacksquare$ \\[0.1cm]
Proposition \ref{prop:geodesicmultwarp} shows that the main difficulty, involved in computing a geodesic $\gamma(t)$ of the multiply-warped Riemannian metric $I$, lies in the solution of the second order differential equation (\ref{eq:geodesic1}). Indeed, once this equation is solved, computing $\gamma(t)$ essentially reduces to an application of $\exp$, which is the Riemannian exponential mapping of the metric $Q$ on $M$. In the context of the present paper, $Q$ is an invariant metric on $M$, where $M$ is a Riemannian symmetric space. Therefore, $\exp$ has a straightforward expression~\cite{helgason} (Theorem 3.3, Page 173). In particular, for the examples treated in Section \ref{sec:examples}, the expression of $\exp$ is well-known in the literature. For the von Mises-Fisher model, this expression is elementary, since geodesics on a sphere in Euclidean space are the great circles on this sphere. For the Riemannian Gaussian model, when this model is defined on the space $M = \mathcal{P}_n$ of $n \times n$ real covariance matrices, the expression of $\exp$ is widely used in the literature, as found in~\cite{pennec2}. \\[0.1cm]
\textbf{Remark 13\,:} the differential equation (\ref{eq:geodesic1}) is the equation of motion of a one-dimensional conservative mechanical system. As such, its solution can be carried out by quadrature~\cite{gallavotti}, (see Page 11). Precisely, the solution reduces to computing the integral
\begin{subequations}
\begin{equation} \label{eq:conservative2}
  t \,=\, \pm\, \int^{\scriptscriptstyle r(t)}_{\scriptscriptstyle r(0)}\frac{dr}{\mathstrut \sqrt{E - V(r)}}
\end{equation}
where the total energy $E$ is a conserved quantity, in the sense that $\dot{E} = 0$, and it can be shown that $E = I_z(U,U)$. Recalling that $dr/d\sigma = \alpha(\sigma)$, this integral can be written
\begin{equation} \label{eq:conservative3}
  t \,=\, \pm\, \int^{\scriptscriptstyle \sigma(t)}_{\scriptscriptstyle \sigma(0)}\frac{\alpha(\sigma)}{\mathstrut \sqrt{E - V(\sigma)}}\,d\sigma \hspace{1cm} 
V(\sigma) = V(r(\sigma))
\end{equation}
Here, if $t$ is interpreted as time, then the integral on the right-hand side gives the time necessary to go from $\sigma(0)$ to $\sigma(t)$. In particular, replacing $\sigma(t)$ by $\infty$ and by $0$ gives the two quantities
\begin{equation} \label{eq:conservative4}
  t_{\scriptscriptstyle \infty} \,=\, \int^{\infty}_{\scriptscriptstyle \sigma(0)}\frac{\alpha(\sigma)}{\mathstrut \sqrt{E - V(\sigma)}}\,d\sigma \hspace{1cm}
 t_{\scriptscriptstyle 0} \,=\, \int^{\scriptscriptstyle \sigma(0)}_{0}\frac{\alpha(\sigma)}{\mathstrut \sqrt{E - V(\sigma)}}\,d\sigma
\end{equation}
where $t_{\scriptscriptstyle \infty}$ is the time necessary for $\sigma(t)$ to reach the value $\sigma = \infty$, and $t_{\scriptscriptstyle 0}$ is the time necessary for $\sigma(t)$ to reach the value $\sigma = 0$. Since $\mathcal{M} \,=\, M\times (0\,,\infty)$, these two values $\sigma = \infty$ and $\sigma = 0$ are excluded from $\mathcal{M}$. Therefore, the geodesic $\gamma(t)$ cannot be extended beyond the time $t = \mathrm{min}(t_{\scriptscriptstyle \infty},t_{\scriptscriptstyle 0})$, as it would then escape from $\mathcal{M}$. \hfill$\blacksquare$ \\[0.1cm]
\end{subequations}
\begin{subequations}
\textbf{Remark 14\,:} a vertical geodesic is a geodesic $\gamma(t)$ for which $\dot{\gamma}(0) = U$ with $U = u_\sigma\,\partial_{\sigma}$. This means that all the $u_q$ are zero. In (\ref{eq:geodesic1}), this implies that $V(r) = 0$, so that $\ddot{r} = 0$ and $r(t)$ is an affine function of $t$. In (\ref{eq:geodesic2}), this implies that $\bar{x}(t) = \bar{x}$ is constant. In Remark 2 of Section \ref{sec:warped}, it was shown that a vertical geodesic is a unique length-minimising geodesic. For a vertical geodesic, (\ref{eq:conservative4}) reads
\begin{equation} \label{eq:conservative5}
t_{\scriptscriptstyle \infty} \,=\, \frac{1}{\mathstrut\sqrt{E}}\,\int^{\infty}_{\scriptscriptstyle \sigma(0)}\alpha(\sigma)\,d\sigma \hspace{1cm}
 t_{\scriptscriptstyle 0} \,=\, \frac{1}{\mathstrut\sqrt{E}}\,\int^{\scriptscriptstyle \sigma(0)}_{0}\alpha(\sigma)\,d\sigma
\end{equation} 
These formulae provide another way of understanding conditions (\ref{eq:complete}) from Section \ref{sec:warped}. Precisely, since $dr/d\sigma = \alpha(\sigma)$, it is clear that $t_{\scriptscriptstyle \infty}$ and $t_{\scriptscriptstyle 0}$ are given by
\begin{equation} \label{eq:conservative6}
t_{\scriptscriptstyle \infty} \,=\, \frac{1}{\mathstrut\sqrt{E}}\, \lim_{\sigma\rightarrow \infty} r(\sigma)- r(\sigma(0)) \hspace{1cm}
t_{\scriptscriptstyle 0} \,=\, \frac{1}{\mathstrut\sqrt{E}}\, \lim_{\sigma\rightarrow 0} r(\sigma(0))- r(\sigma)
\end{equation}
Thus, the first condition in (\ref{eq:complete}) is equivalent to $t_{\scriptscriptstyle \infty} = \infty$, which means that $\sigma(t)$ cannot reach the value $\sigma = \infty$ within a finite time, and the second condition in (\ref{eq:complete}) is equivalent to $t_{\scriptscriptstyle 0} = \infty$, which means that $\sigma(t)$ cannot reach the value $\sigma = 0$ within a finite time. \hfill$\blacksquare$
\end{subequations}
\section{Surprising observation\,: Hadamard manifolds} \label{sec:hadamard}
In Section \ref{sec:warped}, the completeness and curvature of a warped Riemannian metric $I$ were characterised by Formulae (\ref{eq:complete}) and (\ref{eq:curvature}), respectively. Here, based on \ref{subsec:vmf}, these formulae will be applied to the case where $I$ is the Rao-Fisher information metric of the von Mises-Fisher model defined on $S^{n-1}$. Precisely, this application is carried out using a mixture of analytical and numerical computations, for the von Mises-Fisher model defined on $S^{n-1}$ where $n = 2,\ldots,8$. The result is a surprising observation\,: the parameter space of the von Mises-Fisher model, when equipped with the Rao-Fisher information metric $I$, becomes a Hadamard manifold, a simply-connected complete Riemannian manifold of negative sectional curvature~\cite{petersen}\cite{chavel}. Since this observation is true for several values of $n$, it gives rise to a family of Hadamard manifolds. Part of this claim can be proved for any value of $n = 2,\ldots,$ as in the following proposition.  
\begin{proposition} \label{prop:simpcomplete}
  for any value of $n = 2,\ldots,$ the parameter space of the von Mises-Fisher model defined on $S^{n-1}$ is a simply-connected manifold, which moreover becomes a complete Riemannian manifold when equipped with the Rao-Fisher information metric $I$. 
\end{proposition}
\vspace{0.1cm}
\textbf{Proof\,:} recall from Remark 6 that the parameter space of the von Mises-Fisher model defined on $S^{n-1}$ is identified with $\mathbb{R}^n$. Of course, $\mathbb{R}^n$ is a simply-connected manifold~\cite{spanier}. Thus, to prove the proposition, it remains to prove that the parameter space $\mathbb{R}^n$ becomes a complete Riemannian manifold when equipped with the Rao-Fisher information metric $I$. This will be done by proving that all geodesics of the metric $I$ which pass through the point $z = 0$ in $\mathbb{R}^n$ can be extended indefinitely. Then, a corollary of the Hopf-Rinow theorem~\cite{chavel} (Corollary I.7.2, Page 29) implies the required completeness of $\mathbb{R}^n$. 

First, note that the geodesics of the metric $I$ which pass through the point $z = 0$ are exactly the radial straight lines in $\mathbb{R}^n$. Indeed, according to Remark 6, if $\gamma(t)$ is a geodesic of $I$ where $\gamma(t) = (\bar{x}(t),\eta(t))$, then $\gamma(t)$ is identified with the curve $z(t) = \eta(t)\,\bar{x}(t)$ in $\mathbb{R}^n$. Moreover, by Proposition \ref{prop:vmf}, the restriction of $I$ to $\mathbb{R}^n - \lbrace 0 \rbrace$ is a warped Riemannian metric of the general form (\ref{eq:propexp1}). Then, by Remark 14, the vertical geodesics $\gamma(t)$ of this warped Riemannian metric are parameterised by $r(t) =$ affine function of $t$ and $\bar{x}(t) = \bar{x} =$ constant, where $r$ is the vertical distance coordinate. Therefore, each vertical geodesic $\gamma(t)$ can be parameterised by $\eta(t) = \eta(r(t))$ and $\bar{x}(t) = \bar{x} =$ constant. This is identified with the curve $z(t) = \eta(t)\,\bar{x}$, which is a radial straight line in $\mathbb{R}^n$, in the direction $\bar{x}$. It remains to note that the geodesics of the metric $I$ are just the geodesics of its restriction to $\mathbb{R}^n - \lbrace 0 \rbrace$, extended by continuity whenever they reach the point $z = 0$. 

Let $z(t) = \eta(t)\,\bar{x}$ describe a geodesic of the metric $I$, as just explained. To say that this geodesic can be extended indefinitely is equivalent to saying that $\eta(t)$ cannot reach the value $\eta = \infty$ within a finite time. For the von Mises-Fisher model, $\eta(\sigma) = \sigma$ as long as $\sigma \in (0\,,\infty)$. Therefore, according to Remark 14, by evaluating the two conditions (\ref{eq:complete}), it is possible to know whether $\eta(t)$ can reach the two values $\eta = \infty$ and $\eta = 0$ within a finite time. These conditions now read
\begin{subequations}
$$
 \lim_{\eta\rightarrow \infty} r(\eta)- r(\eta_{\scriptscriptstyle 0}) \stackrel{?}{=} \infty \;\;\text{and}\;\;
\lim_{\eta\rightarrow 0} r(\eta_{\scriptscriptstyle 1})- r(\eta) \stackrel{?}{=} \infty
$$
where $\eta_{\scriptscriptstyle 0}$ and $\eta_{\scriptscriptstyle 1}$ are arbitrary. By (\ref{eq:propexp1}), $r(\eta)$ is defined by $dr/d\eta=\left(\,\psi^{\prime\prime}(\eta)\,\right)^{1/2}$. Therefore, the two conditions in (\ref{eq:complete}) are identical to  
\begin{equation} \label{eq:completebis1}
 \int^\infty_{\scriptscriptstyle \eta_{\scriptscriptstyle 0}}\left(\,\psi^{\prime\prime}(\eta)\,\right)^{1/2}\,d\eta \stackrel{?}{=} \infty \;\;\text{and}\;\;
 \int^{\scriptscriptstyle \eta_{\scriptscriptstyle 1}}_{0}\left(\,\psi^{\prime\prime}(\eta)\,\right)^{1/2}\,d\eta \stackrel{?}{=} \infty 
\end{equation} 
where $\psi^{\prime\prime}(\eta)$ is given by (\ref{eq:propvmf1}). For the first integral, recall the asymptotic expansion of modified Bessel functions at $\eta = \infty$~\cite{watson}, (Section 7.23, Page 203. This formula appears with the wrong sign for the second term in parentheses, in~\cite{jupp}),
$$
I_{\nu}(\eta) \,=\, \frac{e^\eta}{\mathstrut \sqrt{\mathstrut 2\pi\eta}}\,\left(\,1 - \frac{4\nu^2-1}{8\eta} + \frac{(4\nu^2-1)(4\nu^2-3^2)}{2!(8\eta)^2}\,\right) \,+\, O\left(\eta^{-3}\,\right)
$$
Using this asymptotic expansion, it follows by performing some direct calculations, and recalling that $\nu = n/2$, 
$$
\begin{array}{ll}
\frac{I_{\nu+1}(\eta)}{\strut I_{\nu+1}(\eta)} \,=& \!\!\!\!1 \,-\, \frac{n}{\eta} \, +\, \frac{n(n-1)}{\strut 2\eta^2} \,+\, O\left(\eta^{-3}\,\right) \\[0.4cm]
\frac{I^2_{\nu}(\eta)}{\strut I^2_{\nu-1}(\eta)} \,= &\!\!\!\! 1 -\, \frac{n-1}{\eta}\,+\, \frac{(n-1)(n-2)}{\strut 2\eta^2} \,+\, O\left(\eta^{-3}\,\right)
\end{array}
$$
Replacing these expressions into (\ref{eq:propvmf1}) immediately gives
\begin{equation} \label{eq:completebis2}
  \psi^{\prime\prime}(\eta) \,=\, \frac{n-1}{2\eta^2} \,+\, O\left(\eta^{-3}\,\right)
\end{equation}
Since $n > 1$, this implies that the first integral in (\ref{eq:completebis1}) is divergent, as required. The second integral in (\ref{eq:completebis1}) is actually convergent. Indeed, $\psi^{\prime\prime}(\eta)$ is a continuous function in the neighborhood of $\eta = 0$, as seen in the proof of Proposition \ref{prop:vmf}, for the limit (\ref{eq:proofvmf4}). Thus, the first condition in (\ref{eq:complete}) is verified, while the second condition is not verified. This means that $\eta(t)$ cannot reach the value $\eta = \infty$ within a finite time, but that it can reach the value $\eta = 0$ within a finite time. The first of these two statements shows that the geodesic described by $z(t)$ can indeed be extended indefinitely. Now, any geodesic of the metric $I$ which passes through the point $z = 0$ is described by some $z(t)$ of this form. \hfill $\blacksquare$\\[0.1cm]
\end{subequations}
\indent The idea behind the proof of the completeness part of Proposition \ref{prop:simpcomplete} can be summarised as follows. The restriction of the Rao-Fisher information metric $I$ to $\mathbb{R}^n - \lbrace 0 \rbrace$ is a warped Riemannian metric. Thus, as stated in Section \ref{sec:warped}, $\mathbb{R}^n - \lbrace 0 \rbrace$ will be a complete Riemannian manifold, when equipped with this warped Riemannian metric, if and only if the two conditions in (\ref{eq:complete}) are verified. Once these conditions are evaluated, it turns out the first one is verified, but the second one is not. Thus, $\mathbb{R}^n - \lbrace 0 \rbrace$ is not a complete Riemannian manifold. However, this is only due to the fact that the point $z = 0$ is excluded. Once this point is included, the parameter space $\mathbb{R}^n$ is obtained, and this is a complete Riemannian manifold, when equipped with the Rao-Fisher information metric $I$. 
Precisely, a vertical geodesic in $\mathbb{R}^n - \lbrace 0 \rbrace$ can reach the point $z = 0$ within a finite time, but then all it does is pass through this point, and immediately return to $\mathbb{R}^n - \lbrace 0 \rbrace$. However, this vertical geodesic cannot escape to infinity within a finite time. With regard to the simple connectedness part of Proposition \ref{prop:simpcomplete}, excluding the point $z = 0$ has no influence if $n > 2$, since $\mathbb{R}^n - \lbrace 0 \rbrace$ is still simply-connected in this case. However, for $n = 2$, $\mathbb{R}^n - \lbrace 0 \rbrace$ is not simply-connected~\cite{spanier}. 

Proposition \ref{prop:simpcomplete} established that the parameter space $\mathbb{R}^n$ of the von Mises-Fisher model is a simply-connected complete Riemannian manifold, for any value of $n$. To show that this parameter space is a Hadamard manifold, it remains to show that it has negative sectional curvature. This is done using numerical computation, for $n = 2,\ldots,8$\,. 

Precisely, let $K_z$ denote the sectional curvature of $\mathbb{R}^n$ at a point $z \in \mathbb{R}^n$, with respect to the Rao-Fisher information metric $I$. Since, according to Proposition \ref{prop:vmf}, $I$ is a warped Riemannian metric, the sectional curvature $K_z$ can be computed from formulae (\ref{eq:curvature}). To evaluate formula (\ref{eq:gauss}), the Gauss equation, it is enough to note that for the von Mises-Fisher model, $K^M$ is a constant equal to $+1$. Indeed, $K^M$ is the sectional curvature of the unit sphere $S^{n-1}$. Then, formula (\ref{eq:gauss}) reads
\begin{subequations} \label{eq:curvaturevmf}
\begin{equation} \label{eq:gaussvmf}
   K_z(u,v) \,=\, \frac{1}{\strut\beta^2}\, -  \left(\frac{\partial_r\beta}{\strut\beta}\right)^2  \hspace{0.2cm} u,v \in T_{\bar{x}}S^{n-1}
\end{equation}
where $z = \eta\,\bar{x}$. On the other hand, formula (\ref{eq:radial}), the radial curvature equation, can be copied directly,
\begin{equation} \label{eq:radialvmf}
K_z(u,\partial_r) \,=\, - \frac{\partial^2_r\beta}{\strut\beta}
\end{equation}
\end{subequations}
Here, $\beta(\eta)$ is given by (\ref{eq:propvmf2}) of Proposition \ref{prop:vmf}, and $r$ is the vertical distance coordinate defined by $dr/d\eta = \left(\,\psi^{\prime\prime}(\eta)\,\right)^{1/2}$. In each one of formulae (\ref{eq:curvaturevmf}), the right-hand side is a real-valued function of $\eta$, independent of the vectors $u, v$. This will be denoted in the following way
\begin{equation} \label{eq:sirfaceradoam}
  K_z(u,v) \,=\, K^s(\eta) \;\;\text{(surface curvature)} \hspace{1cm} K_z(u,\partial_r) \,=\, K^r(\eta)\;\;\text{(radial curvature)} 
\end{equation}
Precisely, $K^s(\eta)$ is the sectional curvature of any section $(u,v)$ tangent to the surface of a sphere centred at $z = 0$ and with Euclidean radius $\eta$. On the other hand, $K^r(\eta)$ is the sectional curvature of a radial section $(u,\partial_r)$. 

Formulae (\ref{eq:curvaturevmf}) were computed numerically for $n = 2,\ldots,8$. It was systematically found that the sectional curvatures $K^s(\eta)$ and $K^r(\eta)$ are negative for all values of $\eta$. From these numerical results, it can be concluded with certainty that the sectional curvature of $\mathbb{R}^n$, with respect to the Rao-Fisher information metric $I$, is negative, when $n$ ranges from $n=2$ to $n=8$. Figures \ref{fig1}--\ref{fig3} below give a graphic representation for $n = 3,4, 5$. The sectional curvatures $K^s(\eta)$ and $K^r(\eta)$ behave in the same way, for all considered values of $n$. Precisely, they are equal to zero at $\eta = 0$, and decrease to a limiting negative value, as $\eta$ becomes large. This limiting value, denoted $K^s(\infty)$ and $K^r(\infty)$, for $K^s(\eta)$ and $K^r(\eta)$, respectively, is given in the following table
\begin{table}[!h]
\begin{center} 
\begin{tabular}{lccccccc}
                    & $n=2$ & $n=3$ & $n=4$ & $n=5$ & $n=6$ & $n=7$ & $n=8$ \\[0.1cm] 
$K^s(\infty)$ &  -0.50 & -0.25 & -0.16 & -0.12 & - 0.10 & -0.08 & -0.07\\[0.12cm]
$K^r(\infty)$ & -0.50 & -0.25 & -0.16 & -0.12 & -0.10 & -0.08 &  -0.07\\
\end{tabular}
\caption{Limiting value of the surface and radial curvatures}
\label{table}
\end{center}
\end{table}
\vspace{-0.4cm}
Remarkably, it appears from this table that $K^s(\infty)$ and $K^r(\infty)$ are close to each other, having the same first two digits after the decimal point. However, it is currently not clear why, or whether it is because $K^s(\infty)$ and $K^r(\infty)$ are equal. 
\vspace{-0.4cm}
\begin{figure} [!h]
\centering
\begin{subfigure}{.5\textwidth}
  \centering
  \includegraphics[width=0.9\linewidth]{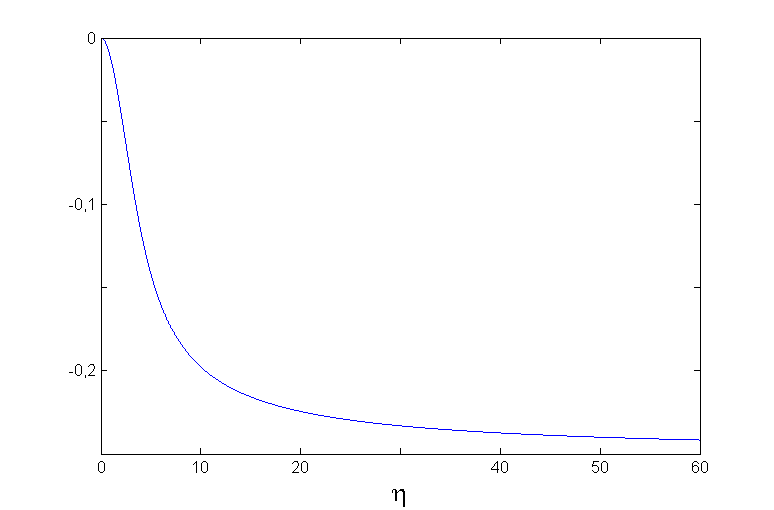}
  \caption{Graph of the surface curvature $K^s(\eta)$}
  \label{fig:hn1}
\end{subfigure}%
\begin{subfigure}{.5\textwidth}
  \centering
  \includegraphics[width=0.9\linewidth]{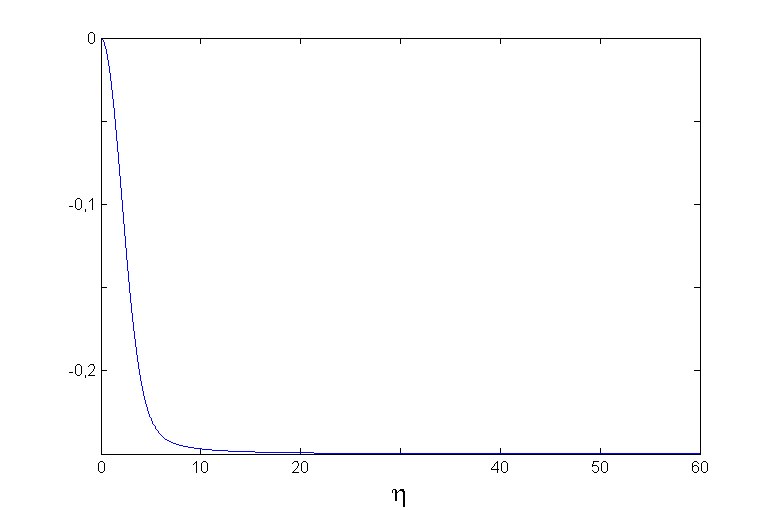}
  \caption{Graph of the radial curvature $K^r(\eta)$}
  \label{fig:hn2}
\end{subfigure}
\caption{($n = 3$) Sectional curvature of the parameter space of the von Mises-Fisher model}
\label{fig1}
\end{figure}
\vspace{-0.05cm}
    
\noindent \textbf{Remark 15\,:} based on Proposition \ref{prop:simpcomplete}, and on the numerical results reported here, it has been found that the parameter space $\mathbb{R}^n$ of the von Mises-Fisher model becomes a Hadamard manifold, when equipped with the Rao-Fisher information metric $I$, for $n = 2,\ldots,8$\,. Indeed, Proposition \ref{prop:simpcomplete} shows that this parameter space is a simply-connected complete Riemannian manifold, for any value of $n = 2,\ldots,$ while the numerical results of Figures \ref{fig1}--\ref{fig3} and Table \ref{table} show that it has negative sectional curvature. Hopefully, future work will soon remove the need to appeal to numerical results, by giving a mathematical proof  of the fact that the sectional curvature of the parameter space $\mathbb{R}^n$ is negative for any value of $n = 2,\ldots,$ without restriction. \hfill$\blacksquare$ \\[0.1cm]
\textbf{Remark 16\,:} preliminary results from ongoing research clearly indicate that the Riemannian Gaussian model, which was studied in \ref{subsec:gauss}, when defined on $M = H^{n-1}$, the $n-1$-dimensional hyperbolic space, has similar properties to the von Mises-Fisher model, with regard to the sectional curvature of its parameter space. Indeed, numerical computations show the sectional curvature of this parameter space $\mathcal{M} = H^{n-1} \times (0\,,\infty)$, equipped with the Rao-Fisher information metric $I$, is negative for $n = 3, 4, 5$. These numerical computations were carried out using formulae (\ref{eq:curvature}), for the sectional curvature of a warped Riemannian metric. This is justified because the hyperbolic space $H^{n-1}$ is an irreducible Riemannian symmetric space, since it is a space of constant negative curvature, so the Rao-Fisher information metric (\ref{eq:proprgm1}) is a warped Riemannian metric. \hfill$\blacksquare$ \\[0.1cm]
\begin{figure} [!t] \label{fig2}
\centering
\begin{subfigure}{.5\textwidth}
  \centering
  \includegraphics[width=0.9\linewidth]{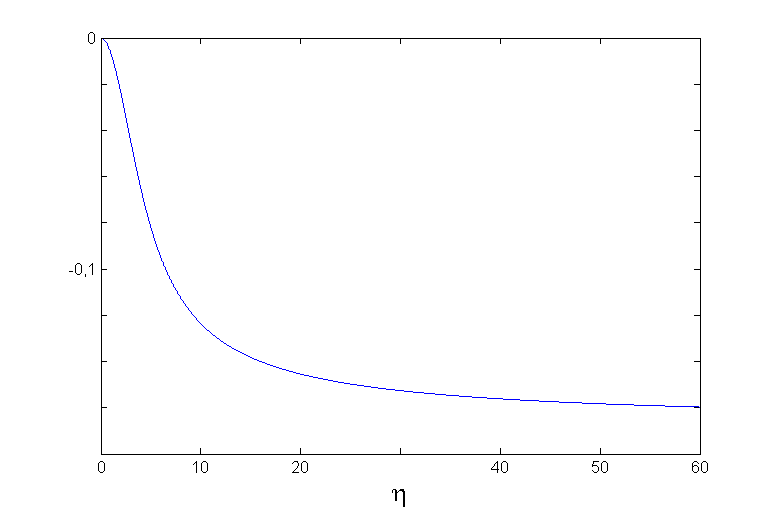}
  \caption{Graph of the surface curvature $K^s(\eta)$}
  \label{fig:hn1}
\end{subfigure}%
\begin{subfigure}{.5\textwidth}
  \centering
  \includegraphics[width=0.9\linewidth]{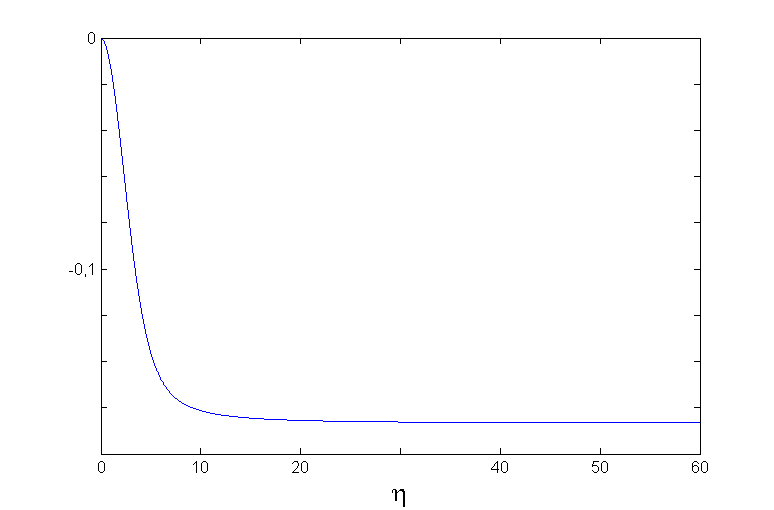}
  \caption{Graph of the radial curvature $K^r(\eta)$}
  \label{fig:hn2}
\end{subfigure}
\caption{($n = 4$) Sectional curvature of the parameter space of the von Mises-Fisher model}
\label{fig2}
\end{figure}
\begin{figure}[!h]
\centering
\begin{subfigure}{.5\textwidth}
  \centering
  \includegraphics[width=0.9\linewidth]{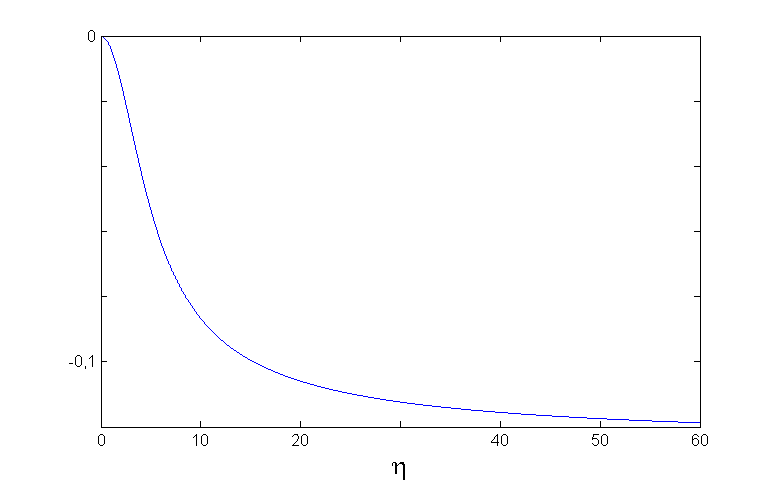}
  \caption{Graph of the surface curvature $K^s(\eta)$}
  \label{fig2}
\end{subfigure}%
\begin{subfigure}{.5\textwidth}
  \centering
  \includegraphics[width=0.9\linewidth]{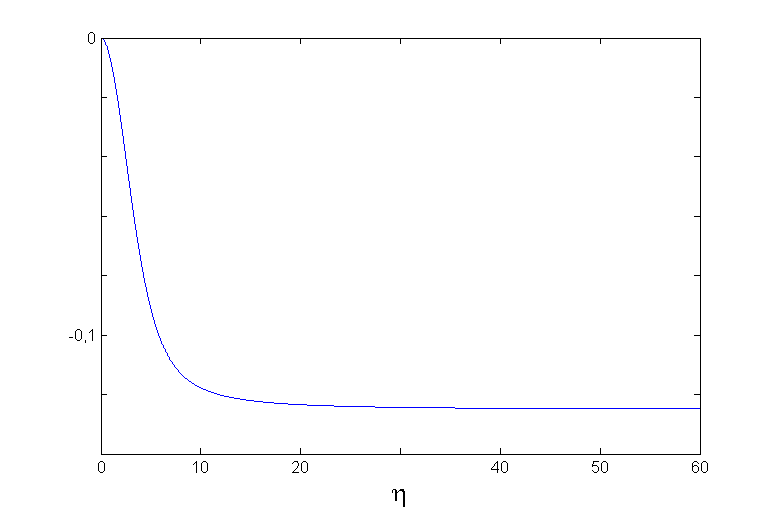}
  \caption{Graph of the radial curvature $K^r(\eta)$}
  \label{fig:hn2}
\end{subfigure}
\caption{($n = 5$) Sectional curvature of the parameter space of the von Mises-Fisher model}
\label{fig3}
\end{figure}

\bibliographystyle{IEEEtran}    
\bibliography{saidetalWARP}

\appendices

\section{Proof of Theorem \ref{th:warplocdis}} \label{app:A}
In order to complete the proof of Theorem \ref{th:warplocdis}, the following proposition is needed. The notation is that of (\ref{eq:fcirc}). 
\begin{proposition} \label{prop:lemma}
\begin{subequations}
  assume condition (\ref{eq:invariance}) holds. Then, 
\begin{equation} \label{eq:lemma1}
  \partial_\sigma\ell(z)\circ g \,=\, \partial_\sigma\ell(g^{-1}\!\cdot z) \hspace{1cm} \nabla_{\bar{x}}\,\ell(z)\circ g \,=\, dg_{\bar{x}}\,\nabla_{\bar{x}}\,\ell(g^{-1}\!\cdot z)\,
\end{equation}
In particular, if $g = s_{\bar{x}}$,
\begin{equation} \label{eq:lemma11}
  \partial_\sigma\ell(z)\circ s_{\bar{x}} \,=\, \partial_\sigma\ell(z) \hspace{1cm} \nabla_{\bar{x}}\,\ell(z)\circ s_{\bar{x}} \,=\, -\,\nabla_{\bar{x}}\,\ell(z)
\end{equation}
\end{subequations}
\end{proposition}
\textbf{Proof\,:} note that (\ref{eq:lemma11}) follows from (\ref{eq:lemma1}),  by the definition of the geodesic-reversing isometry $s_{\bar{x}}$~\cite{helgason}. Indeed, $s_{\bar{x}}\cdot \bar{x} = \bar{x}$ so $s^{-1}_{\bar{x}}\cdot z = z$. Moreover, $ds_{\bar{x}} = -\mathrm{Id}$, as a linear mapping of $T_{\bar{x}}M$, where $\mathrm{Id}$ denotes the identity. To prove (\ref{eq:lemma1}), note that
\begin{subequations}
\begin{equation} \label{eq:prooflem11}
   \left(\partial_\sigma\ell(z)\circ g\right)(x) \,=\, \partial_\sigma\log p(g\cdot x|z) \,=\,
   \partial_\sigma\log p(x|g^{-1}\!\cdot z)
\end{equation}
where the second equality follows from condition (\ref{eq:invariance}). However, 
\begin{equation} \label{eq:prooflem12}
\partial_\sigma\log p(x|g^{-1}\!\cdot z) \,=\, \partial_\sigma\ell(g^{-1}\!\cdot z)(x) 
\end{equation}
Replacing (\ref{eq:prooflem12}) in (\ref{eq:prooflem11}) gives,
$$
\left(\partial_\sigma\ell(z)\circ g\right)(x) \,=\, \partial_\sigma\ell(g^{-1}\!\cdot z)(x)
$$
which is the first part of (\ref{eq:lemma1}).
\end{subequations}
\begin{subequations}
For the second part, a similar reasoning can be applied. Precisely, using condition (\ref{eq:invariance}), it follows, 
\begin{equation} \label{eq:prooflem21}
   \left(d\ell(z)\circ g\right)(x) \,=\, d\log p(g\cdot x|z) \,=\,
   d\log p(x|g^{-1}\!\cdot z) \,=\, d\ell^{(g)}(z)(x)
\end{equation}
where $d\ell(z)$ denotes the derivative of $\ell(z)$ with respect to $\bar{x}$, and $\ell^{(g)}(z) = \ell(g^{-1}\!\cdot z)$, so (\ref{eq:prooflem21}) implies that,
\begin{equation} \label{eq:prooflem22}
   d\ell(z)\circ g \,=\, d\ell^{(g)}(z)
\end{equation}
By the chain rule~\cite{lee}, for $u \in T_{\bar{x}}M$,
$$
\left.d\ell^{(g)}(z)\right|_{\bar{x}}\, u \,=\, d\ell(g^{-1}\!\cdot z) \, dg^{-1}_{\bar{x}}\, u
$$
Replacing in (\ref{eq:prooflem22}),
\begin{equation} \label{eq:prooflem23}
  \left.d\ell(z)\circ g\right|_{\bar{x}} \,=\, d\ell(g^{-1}\!\cdot z) \, dg^{-1}_{\bar{x}}
\end{equation}
\end{subequations}
The second part of (\ref{eq:lemma1}) can now be obtained as follows. By the definition of the Riemannian gradient~\cite{petersen},
\begin{subequations}
\begin{equation} \label{eq:prooflem24}
Q\left(\nabla_{\bar{x}}\,\ell(z)\circ g\,,u\right) \,=\, \left.d\ell(z)\circ g\right|_{\bar{x}}\,u \,=\, d\ell(g^{-1}\!\cdot z) \, dg^{-1}_{\bar{x}}\,u
\end{equation}
where the second equality follows from (\ref{eq:prooflem23}). However, 
$$
d\ell(g^{-1}\!\cdot z) \, dg^{-1}_{\bar{x}}\,u \,=\, Q\left(\nabla_{\bar{x}}\,\ell(g^{-1}\!\cdot z), dg^{-1}_{\bar{x}}\,u\right)
$$
Since $g$ is an isometry of $M$, its derivative $dg_{\bar{x}}$ preserves the Riemannian metric $Q$. Therefore,
\begin{equation} \label{eq:prooflem25}
 Q\left(\nabla_{\bar{x}}\,\ell(g^{-1}\!\cdot z), dg^{-1}_{\bar{x}}\,u\right) \,=\, 
Q\left(dg_{\bar{x}}\,\nabla_{\bar{x}}\,\ell(g^{-1}\!\cdot z),\,u\right)
\end{equation}
Replacing (\ref{eq:prooflem25}) in (\ref{eq:prooflem24}) gives,
$$
Q\left(\nabla_{\bar{x}}\,\ell(z)\circ g\,,u\right) \,=\, Q\left(dg_{\bar{x}}\,\nabla_{\bar{x}}\,\ell(g^{-1}\!\cdot z),\,u\right)
$$
\end{subequations}
To finish the proof, it is enough to note that the vector $u$ is arbitrary. \hfill$\blacksquare$ \\[0.12cm]
\textit{Proof of (\ref{eq:proof2})\,:} recall the polarisation identity, from elementary linear algebra~\cite{lang}, (see Page 29), 
$$
I_z(\partial_\sigma,u) = \frac{1}{4}\,I_z(\partial_\sigma+u,\partial_\sigma+u) - \frac{1}{4}\,I_z(\partial_\sigma-u,\partial_\sigma-u)
$$ 
by replacing (\ref{eq:rao}) into this identity, it can be seen that,  
$$
I_z(\partial_\sigma,u)\,=\,\mathbb{E}_z\left(\,\left(\partial_\sigma\ell(z)\right)\left(d\ell(z)\,u\right)\,\right)
$$
Then, by recalling the definition of the Riemannian gradient~\cite{petersen},
\begin{subequations} 
\begin{equation}
I_z(\partial_\sigma,u)\,=\,\mathbb{E}_z\left(\,\left(\partial_\sigma\ell(z)\right)Q\left(\nabla_{\bar{x}}\,\ell(z),u\right)\,\right)
\end{equation}
Denote the function under the expectation by $f$, and apply (\ref{eq:fcirc}) with $g = s_{\bar{x}}$. Then,
\begin{equation} \label{eq:pa}
\mathbb{E}_z\,f \,=\, \mathbb{E}_{s_{\bar{x}}\cdot z}\,f \,=\, \mathbb{E}_z\left(f\circ s_{\bar{x}}\right)   
\end{equation}
since $s_{\bar{x}}\cdot z = z$. Note that (\ref{eq:proof2}) amounts to saying that $\mathbb{E}_z\,f = 0$. To prove this, note that
\begin{equation}
f\circ s_{\bar{x}} \,=\, \partial_\sigma\ell(z)\circ s_{\bar{x}} \,\times\, Q\left(\nabla_{\bar{x}}\,\ell(z)\circ s_{\bar{x}}\,,u\right) \,=\,
\partial_\sigma\ell(z)\,\times\, - \,Q\left(\nabla_{\bar{x}}\,\ell(z),u\right) \,= -\, f
\end{equation}
\end{subequations}
where the second equality follows from (\ref{eq:lemma11}). Replacing in (\ref{eq:pa}) shows that $\mathbb{E}_z\,f = 0$. \hfill$\blacksquare$ \\[0.1cm]
\textit{Proof of (\ref{eq:proof3})\,:} the idea is to apply Schur's lemma to $I_z(u,u)$, considered as a symmetric bilinear form on $T_{\bar{x}}M$. First, it is shown that this symmetric bilinear form is invariant under the isotropy representation. That is,
\begin{subequations}
\begin{equation} \label{eq:schur1}
  I_z(u,u) \,=\, I_z\left(dk_{\bar{x}}\,u\,,dk_{\bar{x}}\,u\right) \hspace{1cm} \text{for all } k \in K_{\bar{x}}
\end{equation}
This is done using (\ref{eq:fcirc}). Note from (\ref{eq:rao}),
\begin{equation} \label{eq:fb}
I_z(u,u) = \mathbb{E}_z\,\left(\,Q\left(\nabla_{\bar{x}}\,\ell(z),u\right)\,\right)^2
\end{equation}
Denote the function under the expectation by $f$. By (\ref{eq:fcirc}),
\begin{equation} \label{eq:schur2}
\mathbb{E}_z\,f \,=\, \mathbb{E}_{k^{-1}\cdot z}\,f \,=\, \mathbb{E}_z\left(f\circ k^{-1}\right) 
\end{equation}
since $k^{-1}\cdot z = z$ for $k \in K_{\bar{x}}$. To find $f\circ k^{-1}$, note that, 
$$
Q\left(\nabla_{\bar{x}}\,\ell(z)\circ k^{-1},u\right) \,= \,Q\left(dk^{-1}_{\bar{x}}\nabla_{\bar{x}}\,\ell(z),u\right)
\,= \,Q\left(\nabla_{\bar{x}}\,\ell(z),dk_{\bar{x}}\,u\right)
$$
where the first equality follows from (\ref{eq:lemma1}) and the fact that $k\cdot \bar{x} = \bar{x}$, and the second equality from the fact that $dk_{\bar{x}}$ preserves the Riemannian metric $Q$. Now, by (\ref{eq:fb}) and (\ref{eq:schur2}),
$$
I_z(u,u) \,=\, \mathbb{E}_z\,f \,=\, \mathbb{E}_z\left(f\circ k^{-1}\right) \,=\, \mathbb{E}_z\,\left(\,Q\left(\nabla_{\bar{x}}\,\ell(z),dk_{\bar{x}}\,u\right)\,\right)^2 \,=\,  I_z\left(dk_{\bar{x}}\,u\,,dk_{\bar{x}}\,u\right)
$$
and this proves (\ref{eq:schur1}). 
\end{subequations}
\begin{subequations}
Recall Schur's lemma, (\cite{knapp}, Page 240). Applied to (\ref{eq:schur1}), this lemma implies that there exists some multiplicative factor $\beta^2$, such that
\begin{equation} \label{eq:fc}
I_z(u,u) \,=\, \beta^2\,Q_{\bar{x}}(u,u)
\end{equation}
It remains to show that $\beta^2$ is given by (\ref{eq:warplocdis}). Taking the trace of (\ref{eq:fc}),
\begin{equation} \label{eq:fd}
\mathrm{tr}\, I_z = \beta^2\, \mathrm{tr}\,Q_{\bar{x}} = \beta^2\, \mathrm{dim}\,M
\end{equation}
If $e_1\,,\ldots,\,e_d$ is an orthonormal basis of $T_{\bar{x}}M$, then by (\ref{eq:fb}),
\begin{equation} \label{eq:fe}
\mathrm{tr}\, I_z \,= \sum^d_{i=1}\, I_z(e_i,e_i) \,= \mathbb{E}_z\,\sum^d_{i=1}\,  \left(\,Q\left(\nabla_{\bar{x}}\,\ell(z),e_i\right)\,\right)^2 \,= \mathbb{E}_z\,Q\left(\nabla_{\bar{x}}\,\ell(z),\nabla_{\bar{x}}\,\ell(z)\,\right)
\end{equation}
Thus, (\ref{eq:fd}) and (\ref{eq:fe}) show that $\beta^2$ is given by (\ref{eq:warplocdis}).
\end{subequations}
\hfill$\blacksquare$ \\[0.1cm]
To complete the proof of Theorem \ref{th:warplocdis}, it remains to show that the expectations appearing in (\ref{eq:warplocdis}) do not depend on $\bar{x}$. 

For the first expectation, giving $\alpha^2(\sigma)$, note that,
\begin{equation} \label{eq:alpha1}
\mathbb{E}_{g\cdot z}\left(\partial_\sigma\ell(g\cdot z)\right)^2 \,=\, \mathbb{E}_z\left(\partial_\sigma\ell(g\cdot z)\circ g\right)^2 \,=\,
\mathbb{E}_z\left(\partial_\sigma\ell(z)\right)^2
\end{equation}
where the first equality follows  from (\ref{eq:fcirc}) and the second equality follows from (\ref{eq:lemma1}). Thus, this expectation has the same value, whether computed at $g\cdot z = (g\cdot\bar{x},\sigma)$, or at $z = (\bar{x},\sigma)$. Therefore, it does not depend on $\bar{x}$, since the action of $G$ on $M$ is transitive.

For the second expectation, giving $\beta^2(\sigma)$, note that by (\ref{eq:fcirc}),
\begin{subequations}
\begin{equation} \label{eq:beta1}
\mathbb{E}_{g\cdot z} \,Q\left(\nabla_{\bar{x}}\,\ell(g\cdot z)\,,\nabla_{\bar{x}}\,\ell(g\cdot z)\,\right) \,=\, 
\mathbb{E}_{z} \,Q\left(\nabla_{\bar{x}}\,\ell(g\cdot z) \circ g\,,\nabla_{\bar{x}}\,\ell(g\cdot z) \circ g\,\right)
\end{equation}
On the other hand, by (\ref{eq:lemma1}), 
$$
\nabla_{\bar{x}}\,\ell(g\cdot z) \circ g \,=\, dg_{\bar{x}}\,\nabla_{\bar{x}}\,\ell(z) 
$$
Moreover, since $dg_{\bar{x}}$ preserves the Riemannian metric $Q$,
\begin{equation}
Q(\nabla_{\bar{x}}\,\ell(g\cdot z) \circ g\,,\nabla_{\bar{x}}\,\ell(g\cdot z) \circ g) \,=\, Q(dg_{\bar{x}}\,\nabla_{\bar{x}}\,\ell(z)\,,dg_{\bar{x}}\,\nabla_{\bar{x}}\,\ell(z) ) \,= \, Q(\nabla_{\bar{x}}\,\ell(z)\,,\nabla_{\bar{x}}\,\ell(z) )
\end{equation}
Replacing in (\ref{eq:beta1}) gives
\begin{equation} \label{eq:beta2}
 \mathbb{E}_{g\cdot z} \,Q\left(\nabla_{\bar{x}}\,\ell(g\cdot z)\,,\nabla_{\bar{x}}\,\ell(g\cdot z)\,\right) \,=\, \mathbb{E}_z\,Q(\nabla_{\bar{x}}\,\ell(z)\,,\nabla_{\bar{x}}\,\ell(z) )
\end{equation}
\end{subequations}
so this expectation has the same value, at $g\cdot z$ and at $z$. By the same argument made after (\ref{eq:alpha1}), it does not depend on $\bar{x}$. \hfill$\blacksquare$ \\[0.12cm]
\textit{Proof of (\ref{eq:fcirc})\,:} let $dv$ denote the invariant Riemannian volume element of $M$, and note that, 
\begin{subequations}
\begin{equation} \label{eq:ten1}
\mathbb{E}_{g\cdot\, z}\,f \,=\, \int_{M}\, f(x)\,p(x|g\cdot z) \,dv(x) \,=\, \int_{M}\, f(x)\,p(g^{-1}\cdot\!x| z) \,dv(x)
\end{equation}
where the second equality follows from (\ref{eq:invariance}). Introduce the variable $y = g^{-1}\cdot\!x$. Since the volume element $dv$ is invariant,
\begin{equation} \label{eq:ten2}
\int_{M}\, f(x)\,p(g^{-1}\cdot\!x|z) \,dv(x) \,=\, \int_{M}\, f(g\cdot y)\,p(y|z) \,dv(y) \,=\,
\int_{M}\, (f\circ g) (y)\,p(y|z) \,dv(y)
\end{equation}
\end{subequations}
The last integral is the same as $\mathbb{E}_z\left(f\circ g\right)$. Therefore, (\ref{eq:fcirc}) follows from (\ref{eq:ten1}) and (\ref{eq:ten2}). \hfill$\blacksquare$

\section{Proof of Proposition \ref{prop:vmf}} \label{app:B}
It remains to prove (\ref{eq:propvmf1}) and (\ref{eq:propvmf2}). To do so, introduce the following notation, using (\ref{eq:vmf}),
\begin{subequations} \label{eq:prep}
\begin{equation} \label{eq:prep1}
  t = \langle x,\bar{x}\rangle \hspace{0.75cm} Z(\eta) = e^{\psi(\eta)} = (2\pi)^\nu \,\eta^{1-\nu}I_{\nu-1}(\eta)
\end{equation}
Then, $Z(\eta)$ is the moment generating function of $t$, so
\begin{equation} \label{eq:prep2}
  \mathbb{E}_z(t) = \frac{Z^\prime(\eta)}{\strut Z(\eta)} \hspace{0.3cm}\text{and}\hspace{0.3cm} \mathbb{E}_z\left(t^2\right) = \frac{Z^{\prime\prime}(\eta)}{\strut Z(\eta)}
\end{equation}
where the prime denotes differentiation with respect to $\eta$.
\end{subequations}
Recall the derivative and recurrence relations of modified Bessel functions~\cite{watson},
\begin{equation} \label{eq:bessel1}
  \left(\eta^{-a}I_{a}(\eta)\right)^\prime = \eta^{-a}I_{a+1}(\eta) \hspace{0.75cm}  I_{a-1}(\eta) - I_{a+1}(\eta) \,=\, \frac{2a}{\eta}\,I_a(\eta)
\end{equation}
where $a$ is any complex number. By applying these relations to (\ref{eq:prep2}), it is possible to show, through a direct calculation,
\begin{subequations} \label{eq:T}
\begin{eqnarray} 
\label{eq:T1}  \mathbb{E}_z(t) &= \frac{I_\nu(\eta)}{\strut I_{\nu-1}(\eta)}\hspace{1.3cm} \\[0.13cm]
\label{eq:T2} \mathbb{E}_z\left(t^2\right)  & = \frac{1}{n} + \frac{n-1}{n}\, \frac{I_{\nu+1}(\eta)}{\strut I_{\nu-1}(\eta)}
\end{eqnarray}
\end{subequations}
Formulae (\ref{eq:T}) will provide the proof of (\ref{eq:propvmf1}) and (\ref{eq:propvmf2}). \\[0.1cm]
\textit{Proof of (\ref{eq:propvmf1})\,:} since $\psi(\eta)$ is the cumulant generating function of $t$,
\begin{subequations}
\begin{equation} 
\psi^{\prime\prime}(\eta) \,=\, \mathrm{Var}_z(t) \,=\, \mathbb{E}_z\left(t^2\right) - \mathbb{E}_z\left(t\right)^2
\end{equation}
where $\mathrm{Var}$ denotes the variance. Now,  (\ref{eq:propvmf1}) follows immediately by replacing from (\ref{eq:T}) into the right-hand side. \hfill$\blacksquare$ \\[0.12cm]
\textit{Proof of (\ref{eq:propvmf2})\,:} recall from (\ref{eq:proofvmf1}),
\begin{equation} \label{eq:end1}
\beta^2(\eta) \,=\, \frac{\eta^2}{n-1}\,\mathbb{E}_z\left(\, 1 - t^2\,\right)
\end{equation}
However, from (\ref{eq:T2}),
\begin{equation} \label{eq:end2}
\mathbb{E}_z\left(\, 1 - t^2\,\right) = \frac{n-1}{n}\,\left(\, 1 + \frac{I_{\nu+1}(\eta)}{\strut I_{\nu-1}(\eta)}\,\right)
\end{equation}
\end{subequations}
Now, (\ref{eq:propvmf2}) follows by replacing (\ref{eq:end2}) into (\ref{eq:end1}). \hfill $\blacksquare$ \\[0.12cm]
\textit{Proof of (\ref{eq:T1})\,:} using the derivative relation of modified Bessel functions, which is the first relation in (\ref{eq:bessel1}), with $a = \nu - 1$, it follows that
\begin{subequations}
\begin{equation} \label{eq:z1}
Z^\prime(\eta) \,=\, (2\pi)^\nu\,\eta^{1-\nu}I_\nu(\eta)
\end{equation}
Now, Formula (\ref{eq:T1}) follows by replacing this into (\ref{eq:prep2}) and using (\ref{eq:prep1}). \hfill $\blacksquare$ \\[0.1cm]
\textit{Proof of (\ref{eq:T2})\,:} write (\ref{eq:z1}) in the form
$$
Z^\prime(\eta) \,=\, (2\pi)^\nu\,\eta\,\left(\eta^{-\nu}I_\nu(\eta)\right)
$$
By the product rule
$$
Z^{\prime\prime}(\eta) \,=\, (2\pi)^\nu\, \eta^{-\nu}I_\nu(\eta) + (2\pi)^\nu\,\eta\,\left(\eta^{-\nu}I_\nu(\eta)\right)^\prime
$$
The derivative in the second term can be evaluated from the derivative relation of modified Bessel functions, with $a = \nu$. Then,
$$
Z^{\prime\prime}(\eta) \,=\, (2\pi)^\nu\, \eta^{-\nu}I_\nu(\eta) + (2\pi)^\nu\,\eta^{1-\nu}I_{\nu+1}(\eta)
$$
Rearrange this formula as
$$
Z^{\prime\prime}(\eta) \,=\, (2\pi)^\nu\,\eta^{1-\nu}\,\left( \eta^{-1}I_\nu(\eta) + I_{\nu+1}(\eta)\right)
$$
By the recurrence relation of modified Bessel functions, which is the second relation in (\ref{eq:bessel1}), with $a = \nu$, it then follows
$$
Z^{\prime\prime}(\eta) \,=\, (2\pi)^\nu\,\eta^{1-\nu}\,\left( \frac{1}{2\nu}\,I_{\nu-1}(\eta) - \frac{1}{2\nu}\,I_{\nu+1}(\eta) + I_{\nu+1}(\eta)\right)
$$
Recalling that $2\nu = n$, this can be written,
\begin{equation} \label{eq:z2}
 Z^{\prime\prime}(\eta) \,=\, (2\pi)^\nu\,\eta^{1-\nu}\,\left( \frac{1}{n}\,I_{\nu-1}(\eta) + \frac{n-1}{n}\, I_{\nu+1}(\eta)\right)
\end{equation}
Now, Formula (\ref{eq:T2}) follows by replacing this into (\ref{eq:prep2}) and using (\ref{eq:prep1}). \hfill $\blacksquare$ \\[0.1cm]

\end{subequations}

\section{Proof of Proposition \ref{prop:geodesicmultwarp}} \label{app:C}
The setting and notations are the same as in Section \ref{sec:geodesic}, except for the fact that $\bar{x}$ is written as $x$, without the bar, in order to avoid notations such as $\dot{\bar{x}}$ or $\ddot{\bar{x}}$. This being said, let $\tilde{\nabla}$ and $\nabla$ denote the Levi-Civita connections of the Riemannian metrics $I$ and $Q$, respectively. Thus, $\tilde{\nabla}$ is a connection on the tangent bundle of the manifold $\mathcal{M}$, and $\nabla$ is a connection on the tangent bundle of the manifold $M$\!~\cite{petersen}\cite{chavel}. Introduce the shape operator $S:T_{x}M\rightarrow T_{x}M$, which is given as in~\cite{petersen},
\begin{equation} \label{eq:shape}
  S(u) = \tilde{\nabla}_u\,\partial_r \hspace{1cm} u \in T_{x}M
\end{equation}
for any $x \in M$. The following identities can be found in~\cite{petersen} (Section 2.4, Page 41),
\begin{subequations} \label{eq:shapeid}
\begin{eqnarray}
\label{eq:shape1}\hspace{1.7cm}&  \tilde{\nabla}_{\partial_r}\,\partial_r \!\!\!\!\!& = 0\hspace{5cm} \\[0.1cm]
\label{eq:shape2}& \tilde{\nabla}_{\partial_r}\,X \!\!\!\!\!&=  S(X)\hspace{4.4cm} \\[0.1cm]
\label{eq:shape3}&\!\! \tilde{\nabla}_{X}\,Y \!\!\!\!\!&= \nabla_X\,Y - I(S(X),Y)\,\partial_r \hspace{1.9cm}\,
\end{eqnarray}
for any vector fields $X$ and $Y$ on $M$. Using these identities, it is possible to write the geodesic equation of the Riemannian metric $I$, in terms of the shape operator $S$. This is given in the following proposition.
\vspace{0.1cm}

\end{subequations}
\begin{proposition} \label{prop:geodesicequation1}
 let $\gamma(t)$ be a curve in $\mathcal{M}$\,, with $\gamma(t) = (x(t),\sigma(t))$ and let $r(t) = r(\sigma(t))$. The curve $\gamma(t)$ is a geodesic of the Riemannian metric $I$ if and only if it satisfies the geodesic equation
\begin{subequations} \label{eq:geoeq} 
 \begin{eqnarray}
 \ddot{r} \!\! &=  I(S(\dot{x}),\dot{x}) \\[0.1cm]
 \,\ddot{x} \!\! &= -2\,\dot{r}S(\dot{x}) \hspace{0.1cm}
 \end{eqnarray}
\end{subequations}
where $\ddot{x} = \nabla_{\dot{x}}\,\dot{x}$ is the acceleration of the curve $x(t)$ in $M$.
\end{proposition}
\vspace{0.1cm}
The shape operator $S$ moreover admits a simple expression, which can be derived from expression (\ref{eq:mwarpedbis}) of the Riemannian metric $I$, using the fact that the Levi-Civita connection $\tilde{\nabla}$ is a metric connection~\cite{chavel} (Theorem I.5.1, Page 16).
\vspace{0.1cm}

\begin{proposition} \label{prop:shape}
in the notation of (\ref{eq:mwarpedbis}), the shape operator $S$ is given by
\begin{equation} \label{eq:propshape}
  S(u) = \sum^r_{q=1}\, \frac{\partial_r\,\beta_q(r)}{\strut \beta_q(r)}\,u_q
\end{equation} 
In other words, the decomposition $u = u_1\,+\ldots+\,u_r$ provides a block-diagonalisation of $S$, where each block is a multiple of identity. 
\end{proposition}
\vspace{0.1cm}
Combining Propositions \ref{prop:geodesicequation1} and \ref{prop:shape}, the geodesic equation (\ref{eq:geoeq}) takes on a new form. Precisely, replacing (\ref{eq:propshape}) into (\ref{eq:geoeq}) gives the following equations
\begin{subequations} \label{eq:geoeq1} 
 \begin{eqnarray}
\label{eq:geoeq11}& \!\!\!\ddot{r} \!\!\!\! &=  \sum^r_{q=1}\, \beta_q(r)\partial_r\beta_q(r)\,Q(\dot{x}_q,\dot{x}_q) \\[0.12cm]
\label{eq:geoeq21}& \ddot{x}_q  \!\!\!\! &= -2\,\dot{r}\,\frac{\partial_r\,\beta_q(r)}{\mathstrut \beta_q(r)}\,\dot{x}_q \hspace{0.1cm} 
 \end{eqnarray}
\end{subequations}
where $\dot{x} = \dot{x}_1\,+\ldots+\,\dot{x}_r$ and $\ddot{x} = \ddot{x}_1\,+\ldots+\,\ddot{x}_{r\,}$. The proof of Proposition \ref{prop:geodesicmultwarp} can be obtained directly from equations (\ref{eq:geoeq1}), using the following conservation laws.
\vspace{0.1cm}

\begin{proposition} \label{prop:conservation}
  each one of the following quantities $C_q$ is a conserved quantity, 
\begin{equation} \label{eq:conservation11}
  C_q \,=\, \beta^4_q(r)\,Q(\dot{x}_{q\,},\dot{x}_q)\hspace{1cm} \text{for }\,\, q = 1\,,\ldots,\,r
\end{equation} 
In other words, $C_q$ remains constant when evaluated along any geodesic $\gamma(t)$ of the Riemannian metric $I$. 
\end{proposition}
\vspace{0.1cm}

For now, assume that Propositions \ref{prop:geodesicequation1} through \ref{prop:conservation} are true. To prove Proposition \ref{prop:geodesicmultwarp}, note the following.\\[0.1cm]
\textit{Proof of (\ref{eq:geodesic1})\,:} it is enough to show that the right-hand side of (\ref{eq:geodesic1}) is the same as the right-hand side of (\ref{eq:geoeq11}). To do so, note from (\ref{eq:geodesic1}) and (\ref{eq:conservation11}) that
\begin{subequations}
\begin{equation} \label{eq:potentials}
V(r) \,=\, \sum^r_{q=1}\,\frac{\beta^2_q(r(0))}{\beta^2_q(r)}\,I_z(u_{q\,},u_q) \,=\, \sum^r_{q=1}\,\frac{C_q}{\strut \beta^2_q(r)}
\end{equation}
Indeed, since $\dot{x}_q(0) = u_q$ and since $C_q$ is a conserved quantity
$$
\beta^2_q(r(0))\,I_z(u_{q\,},u_q) \,=\, \left.\beta^4_q(r)\,Q(\dot{x}_{q\,},\dot{x}_q)\right|_{t=0} \,=\,C_q 
$$
Now, replacing the derivative of (\ref{eq:potentials}) into the right-hand side of (\ref{eq:geodesic1}) directly leads to the right-hand side of (\ref{eq:geoeq11}). \hfill$\blacksquare$\\[0.1cm]
\textit{Proof of (\ref{eq:geodesic2})\,:} recall from Remark 7 that $M$ is the Riemannian product of the $M_q$. Therefore, the Riemannian exponential mapping of $M$ is also the product of the Riemannian exponential mappings of the $M_q$. Precisely, (\ref{eq:geodesic2}) is equivalent to
\begin{equation} \label{eq:reparam}
x_q(t) \,=\, \exp_{\scriptscriptstyle {x_q(0)}}\,\left[\,\left(\,\int^t_0 \frac{\beta^2_q(r(0))}{\strut \beta^2_q(r(s))}ds\right)\,u_q\,\right] \hspace{1cm} \text{for }\,\, q = 1\,,\ldots,\,r
\end{equation}
This means that the curve $x_q(t)$ in $M_q$ is a reparameterised geodesic $\left(\delta_q\circ F\right)(t)$ where $\delta_q(t)$ is the geodesic given by $\delta_q(t)= \exp(t\,u_q)$ and $F(t)$ is the integral inside the parentheses in (\ref{eq:reparam}). To prove (\ref{eq:geodesic2}), it is sufficient to prove that (\ref{eq:reparam}) solves equation (\ref{eq:geoeq21}). Using the chain rule,  (\ref{eq:reparam}) implies that
\vfill
\pagebreak
\begin{equation}\label{eq:rereparam}
\ddot{x}_q \,=\, \dot{F}^{\,2}\left(\ddot{\delta}_q\circ F\right)\,+\, F^{\,\prime\prime}\left(\dot{\delta}_q\circ F\right) \,=\, \dot{F}^{\,2}\left(\ddot{\delta}_q\circ F\right)\,+\, \frac{F^{\,\prime\prime}}{F^\prime}\,\dot{x}_q \,=\, \frac{F^{\,\prime\prime}}{F^\prime}\,\dot{x}_q 
\end{equation}
where the third equality follows because $\delta_q$ is a geodesic, and therefore its acceleration $\ddot{\delta}_q$ is zero. By replacing the definition of the function $F(t)$, it is seen that (\ref{eq:rereparam}) is the same as (\ref{eq:geoeq21}). It follows that (\ref{eq:reparam}) solves (\ref{eq:geoeq21}), as required.
\end{subequations}
\hfill$\blacksquare$ \\[0.1cm]
The proofs of propositions \ref{prop:geodesicequation1} through \ref{prop:conservation} are now given. \\[0.1cm]
\textbf{Proof of Proposition \ref{prop:geodesicequation1}\,:} recall the geodesic equation is $\tilde{\nabla}_{\dot{\gamma}}\,\dot{\gamma} \,=0$, which means that the velocity $\dot{\gamma}(t)$ is self-parallel\!~\cite{petersen}\cite{chavel}. Here, the velocity $\dot{\gamma}(t)$ is given by $\dot{\gamma}_t \,=\, \dot{r}\,\partial_r \,+\, \dot{x}$. Accordingly, the left-hand side of the geodesic equation is
\begin{subequations}
\begin{equation} \label{eq:prop91}
\tilde{\nabla}_{\dot{\gamma}}\,\dot{\gamma} \,=\, \tilde{\nabla}_{\dot{\gamma}}\,\dot{r}\,\partial_r  \,+\, \tilde{\nabla}_{\dot{\gamma}}\,\dot{x} \,=\, \ddot{r}\,\partial_r \,+\,\dot{r}\,\tilde{\nabla}_{\dot{\gamma}}\,\partial_r\,+\, \tilde{\nabla}_{\dot{\gamma}}\,\dot{x}  
\end{equation}
where the second equality follows by the product rule for the covariant derivative\!~\cite{petersen}\cite{chavel}. The second and third terms on the right-hand side of (\ref{eq:prop91}) can be written in terms of the shape operator $S$. Precisely, for the second term,
\begin{equation} \label{eq:prop92}
  \tilde{\nabla}_{\dot{\gamma}}\,\partial_r \,=\, \dot{r}\,\tilde{\nabla}_{\partial_r}\,\partial_r\,+\, \tilde{\nabla}_{\dot{x}}\,\partial_r \,=\, S(\dot{x})
\end{equation}
where the second equality follows from (\ref{eq:shape}) and (\ref{eq:shape1}). Moreover, for the third term,
\begin{equation} \label{eq:prop93}
  \tilde{\nabla}_{\dot{\gamma}}\,\dot{x} \,=\, \dot{r}\,\tilde{\nabla}_{\partial_r}\,\dot{x}\,+\, \tilde{\nabla}_{\dot{x}}\,\dot{x} \,=\, \dot{r}\,S(\dot{x}) \,+\,\ddot{x} \,-\, I(S(\dot{x}),\dot{x})\,\partial_r
\end{equation}
where the second equality follows from (\ref{eq:shape2}) and (\ref{eq:shape3}). Replacing (\ref{eq:prop92}) and (\ref{eq:prop93}) into (\ref{eq:prop91}), the left-hand side of the geodesic equation becomes
$$
\tilde{\nabla}_{\dot{\gamma}}\,\dot{\gamma} \,=\, \left(\,\ddot{r} - I(S(\dot{x}),\dot{x})  \,\right)\,\partial_r \,+\, \left(\, \ddot{x} + 2\dot{r}S(\dot{x})\,\right)
$$
Setting this equal to zero immediately gives equations (\ref{eq:geoeq}). 
\end{subequations}
\hfill$\blacksquare$\\[0.1cm]
\textbf{Proof of Proposition \ref{prop:shape}\,:} recall the shape operator $S$ is symmetric, since it is essentially the Riemannian Hessian of $r$~\cite{petersen} (Section 2.4, Page 41). Therefore, it is enough to evaluate $I(S(u),u)$ for $u \in T_xM$. Let $X$ by a vector field on $M$, with $X(x) = u$. Then,
\begin{subequations}
\begin{equation} \label{eq:i1}
I(S(u),u) \,=\, I\left(S(X),X\right)\,=\,I\left(\,\tilde{\nabla}_{\partial_r}\,X\,,X\right)
\end{equation}
where the second equality follows from (\ref{eq:shape2}). Using the fact that the Levi-Civita connection $\tilde{\nabla}$ is a metric connection\!~\cite{chavel} (Theorem I.5.1, Page 16), the right-hand side can be written as
\begin{equation} \label{eq:i2}
I\left(\,\tilde{\nabla}_{\partial_r}\,X\,,X\right) \,=\, \frac{1}{2}\, \partial_r\,I(X,X) \,=\, \frac{1}{2}\, \partial_r\,\sum^r_{q=1}\,\beta^2_q(r)\,Q_x(u_{q\,},u_q) 
\end{equation}
where the second equality follows from (\ref{eq:mwarpedbis}). It remains to note that
$$
\frac{1}{2}\, \partial_r\,\beta^2_q(r)\,Q_x(u_{q\,},u_q) \,=\, \frac{\partial_r\,\beta_q(r)}{\strut \beta_q(r)}\,I_z(u_{q\,},u_q)
$$
Accordingly, (\ref{eq:i1}) and (\ref{eq:i2}) imply
\begin{equation} \label{eq:i3}
I(S(u),u) \,=\, I\left(\, \sum^r_{q=1}\, \frac{\partial_r\,\beta_q(r)}{\strut \beta_q(r)}\,u_{q\,},u_q \,\right) \,=\,
I\left(\, \sum^r_{q=1}\, \frac{\partial_r\,\beta_q(r)}{\strut \beta_q(r)}\,u_{q\,},u \,\right)
\end{equation}
and (\ref{eq:propshape}) follows from the fact that $S$ is symmetric. 
\end{subequations}
\hfill$\blacksquare$\\[0.1cm]
\textbf{Proof of Proposition \ref{prop:conservation}\,:} to say that $C_q$ is a conserved quantity means that $\dot{C}_q = 0$. From (\ref{eq:conservation11}),
\begin{subequations}
\begin{equation} \label{eq:provcons1}
\dot{C}_q \,=\, 4\dot{r}\,\beta^3_q(r)\,\partial_r\beta_q(r)\,Q(\dot{x}_{q\,}\,\dot{x}_{q}) \,+\,\beta^4_q(r)\,\frac{d}{dt}\,Q(\dot{x}_{q\,}\,\dot{x}_{q})
\end{equation}
The last derivative can be expressed as
\begin{equation} \label{eq:provcons2}
\frac{d}{dt}\,Q(\dot{x}_{q\,}\,\dot{x}_{q}) \,=\, 2\,Q(\ddot{x}_{q\,},\dot{x}_q) \,=\, -4\dot{r}\,\frac{\partial_r\beta_q(r)}{\strut \beta_q(r)}\,Q(\dot{x}_{q\,},\dot{x}_q)
\end{equation} 
where the second equality follows from (\ref{eq:geoeq21}). By replacing (\ref{eq:provcons2}) into (\ref{eq:provcons1}), it follows immediately that $\dot{C}_q = 0$.
\end{subequations}
\hfill$\blacksquare$



\end{document}